\documentclass[12pt]{article}
\usepackage{amsmath,amsfonts,epsfig}
\title{\bf Matching Theorems for Systems of a Finitely Generated Coxeter Group}
\author{Michael Mihalik, John Ratcliffe, and Steven Tschantz \\ 
Mathematics Department, Vanderbilt University, \\
Nashville TN 37240, USA}
\newtheorem{theorem}{Theorem}[section]
\newtheorem{proposition}[theorem]{Proposition}
\newtheorem{lemma}[theorem]{Lemma}

\newenvironment{example}{\vspace{.2in}{\noindent\bf Example\ }}{\vspace{.2in}}
\newenvironment{proof}{{\bf Proof:\ }}{\hfill$\square$\vspace{.2in}}

\newcommand{\realnos}{{\mathbb R}}

\def\diag{{\rm diag}}
\def\ov{\overline}
\def\lra{\longrightarrow}
\date{}
\begin{document}
\maketitle

\section{Introduction} 

The isomorphism problem for finitely generated Coxeter groups is the problem of deciding 
if two finite Coxeter matrices define 
isomorphic Coxeter groups. Coxeter \cite{Coxeter} solved this problem for finite irreducible 
Coxeter groups.  
Recently there has been considerable interest and activity on the isomorphism problem 
for arbitrary finitely generated Coxeter groups. 
For a recent survey, see M\"uhlherr \cite{Muhlherr}. 

The isomorphism problem for finitely generated Coxeter groups 
is equivalent to the problem of determining all 
the automorphism equivalence classes of 
sets of Coxeter generators for an arbitrary finitely generated Coxeter group. 
In this paper,  we prove a series of matching theorems for 
two sets of Coxeter generators of a finitely generated Coxeter group 
that identify common features of the two sets of generators.  
As an application, we describe an algorithm for finding a set of 
Coxeter generators of maximum rank 
for a finitely generated Coxeter group.  

In \S 2, we state some preliminary results. 
In \S 3, we prove a matching theorem for two systems of a finite Coxeter group. 
In \S 4, we prove our Basic Matching Theorem between the sets of maximal noncyclic irreducible 
spherical subgroups of two systems of a finitely generated Coxeter group. 
In \S 5, we study nonisomorphic basic matching. 
In \S 6, we prove a matching theorem between the sets of noncyclic irreducible spherical subgroups 
of two systems of a finitely generated Coxeter group.  
As an application, we prove the Edge Matching Theorem. 
In \S 7, we discuss twisting and flattening visual graph of groups decompositions of Coxeter systems. 
In \S 8, we prove the Decomposition Matching Theorem. 
In \S 9, we prove the Simplex Matching Theorem. 
In \S 10, we describe our algorithm for finding a set of 
Coxeter generators of maximum rank 
for a finitely generated Coxeter group.  

\section{Preliminaries} 
 
A {\it Coxeter matrix} is a symmetric matrix $M = (m(s,t))_{s,t\in S}$ 
with $m(s,t)$ either a positive integer or infinity and $m(s,t) =1$ 
if and only if $s=t$. A {\it Coxeter system} with Coxeter matrix $M = (m(s,t))_{s,t\in S}$ 
is a pair $(W,S)$ consisting of a group $W$ and a set of generators $S$ for $W$ 
such that $W$ has the presentation
$$W =\langle S \ |\ (st)^{m(s,t)}:\, s,t \in S\ \hbox{and}\ m(s,t)<\infty\rangle$$
If $(W,S)$ is a Coxeter system with Coxeter matrix $M = (m(s,t))_{s,t\in S}$, 
then the order of $st$ is $m(s,t)$ for each $s,t$ in $S$ by 
Prop. 4, p. 92 of Bourbaki \cite{Bourbaki}, and so a Coxeter system $(W,S)$ 
determines its Coxeter matrix; moreover, any Coxeter matrix 
$M = (m(s,t))_{s,t\in S}$ determines a Coxeter system $(W,S)$ where $W$ 
is defined by the above presentation. If $(W,S)$ is a Coxeter system, 
then $W$ is called a {\it Coxeter group} and $S$ is called a set of {\it Coxeter generators} 
for $W$, and the cardinality of $S$ is called the {\it rank} of $(W,S)$.

\begin{proposition} 
A Coxeter system $(W,S)$ has finite rank if and only if $W$ 
is finitely generated. 
\end{proposition}
\begin{proof}
This follows Theorem 2 (iii), p. 20 of Bourbaki  \cite{Bourbaki}. 
\end{proof}

Let $(W,S)$ be a Coxeter system. A {\it visual subgroup} of $(W,S)$ 
is a subgroup of $W$ of the form $\langle A\rangle$ for some $A \subset S$. 
If $\langle A\rangle$ is a visual subgroup of $(W,S)$, 
then $(\langle A\rangle, A)$ is also a Coxeter system 
by Theorem 2 (i), p. 20 of Bourbaki  \cite{Bourbaki}.

When studying a Coxeter system $(W,S)$ with Coxeter matrix $M$ 
it is helpful to have a visual representation of $(W,S)$. 
There are two graphical ways of representing $(W,S)$ 
and we shall use both depending on our needs. 

The {\it Coxeter diagram} ({\it {\rm C}-diagram}) of $(W,S)$ is the labeled undirected graph 
$\Delta = \Delta(W,S)$ with vertices $S$ and edges 
$$\{(s,t) : s, t \in S\ \hbox{and}\ m(s,t) > 2\}$$
such that an edge $(s,t)$ is labeled by $m(s,t)$. 
Coxeter diagrams are useful for visually representing finite Coxeter groups. 
If $A\subset S$, then $\Delta(\langle A\rangle,A)$ is the 
subdiagram of $\Delta(W,S)$ induced by $A$.

A Coxeter system $(W,S)$ is said to be {\it irreducible} 
if its C-diagram $\Delta$ is connected. 
A visual subgroup $\langle A\rangle$ of $(W,S)$ is said to be {\it irreducible} 
if $(\langle A\rangle, A)$ is irreducible. 
A subset $A$ of $S$ is said to be {\it irreducible} if $\langle A\rangle$ is irreducible. 

A subset $A$ of $S$ is said to be a {\it component} of $S$ if $A$ is a maximal irreducible 
subset of $S$ or equivalently if $\Delta(\langle A\rangle, A)$ is a connected component 
of $\Delta(W,S)$. 
The connected components of the $\Delta(W,S)$
represent the factors of a direct product decomposition of $W$.

The {\it presentation diagram} ({\it {\rm P}-diagram}) of $(W,S)$ is the labeled undirected graph 
$\Gamma = \Gamma(W,S)$ with vertices $S$ and edges 
$$\{(s,t) : s, t \in S\ \hbox{and}\ m(s,t) < \infty\}$$
such that an edge $(s,t)$ is labeled by $m(s,t)$. 
Presentation diagrams are useful for visually representing infinite Coxeter groups. 
If $A\subset S$, then $\Gamma(\langle A\rangle,A)$ is the 
subdiagram of $\Gamma(W,S)$ induced by $A$. 
The connected components of $\Gamma(W,S)$
represent the factors of a free product decomposition of $W$.

\begin{example} 
Consider the Coxeter group $W$ generated by the four reflections 
in the sides of a rectangle in $E^2$.  The C-diagram of $(W,S)$ is 
the disjoint union of two edges labeled by $\infty$.  
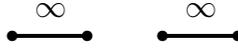
$$\mbox{
\setlength{\unitlength}{1cm}
\begin{picture}(3,.5)(0,0)
\thicklines
\put(0,0){\line(1,0){1}}
\put(2,0){\line(1,0){1}}
\put(0,0){\circle*{.15}}
\put(1,0){\circle*{.15}}
\put(2,0){\circle*{.15}}
\put(3,0){\circle*{.15}}
\put(.3,.25){$\infty$}
\put(2.3,.25){$\infty$}
\end{picture}}$$
Therefore $W$ is the direct product of two infinite dihedral groups.  
The P-diagram of $W$ is a square with edge labels 2.   
$$\mbox{
\setlength{\unitlength}{1cm}
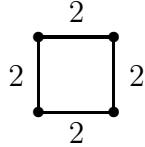
\begin{picture}(1.5,1.5)(0,0)
\thicklines
\put(0,0){\line(0,1){1}}
\put(0,0){\line(1,0){1}}
\put(1,0){\line(0,1){1}}
\put(0,1){\line(1,0){1}}
\put(0,0){\circle*{.15}}
\put(1,0){\circle*{.15}}
\put(1,1){\circle*{.15}}
\put(0,1){\circle*{.15}}
\put(-.4,.35){2}
\put(.4,1.2){2}
\put(1.2,.35){2}
\put(.4,-.4){2}
\end{picture}}$$
\end{example} 

Let $(W,S)$ and $(W',S')$ be Coxeter systems 
with P-diagrams $\Gamma$ and $\Gamma'$, respectively. 
An {\it isomorphism} $\phi: (W,S) \to (W',S')$ of Coxeter systems 
is an isomorphism $\phi: W\to W'$ such that $\phi(S) = S'$. 
An {\it isomorphism} $\psi:\Gamma\to \Gamma'$ of P-diagrams is a bijection 
from $S$ to $S'$ that preserves edges and their labels. 
\begin{proposition} 
Let $(W,S)$ and $(W',S')$ be Coxeter systems with P-dia\-grams $\Gamma$ and $\Gamma'$, 
respectively. Then
\begin{enumerate}
\item $(W,S) \cong (W',S')$ if and only if $\Gamma \cong \Gamma'$,
\item $W\cong W'$ if and only if $W$ has a set of Coxeter generators $S''$ such that 
		$(W,S'')\cong (W',S')$,
\item $W\cong W'$ if and only if $W$ has a P-diagram $\Gamma''$ such that 
		$\Gamma''\cong \Gamma'$.
\end{enumerate}
\end{proposition} 
\begin{proof} (1) If $\phi: (W,S) \to (W',S')$ is an isomorphism, then 
$\phi$ restricts to an isomorphism $\overline\phi:\Gamma \to \Gamma'$ and 
if $\psi:\Gamma \to \Gamma'$ is an isomorphism, then 
$\psi$ extends to a unique isomorphism $\hat \psi: (W,S) \to (W',S')$. 

(2) If $\phi:W\to W'$ is an isomorphism, then $S''=\phi^{-1}(S')$ 
is a set of Coxeter generators for $W$ and $\phi:(W,S'')\to (W',S')$ 
is an isomorphism. 

(3) Statement (3) follows from (1) and (2).
\end{proof}
\begin{proposition} 
Let $W$ be a Coxeter group and let ${\cal S}$ be the collection 
of sets of Coxeter generators for $W$. Then
\begin{enumerate}
\item The group ${\rm Aut}(W)$ acts on ${\cal S}$.
\item Sets of Coxeter generators $S$ and $S'$ for $W$ are in the same ${\rm Aut}(W)$-orbit 
		if and only if $(W,S) \cong (W,S')$. 
\item The set of ${\rm Aut}(W)$-orbits ${\cal S}/{\rm Aut}(W)$ is in one-to-one 
		correspondence with the set of isomorphism classes of P-diagrams for $W$. 
\end{enumerate}
\end{proposition}
\begin{proof} (1) and (2) are obvious. (3) follows from (2) and Prop. 2.2(1).
\end{proof}

A Coxeter group $W$ is said to be {\it rigid} 
if for any two sets of Coxeter generators $S$ and $S'$ for $W$, 
there is an automorphism $\alpha: (W,S)\to (W,S')$ 
or equivalently any two sets of Coxeter generators $S$ and $S'$ for $W$ 
determine isomorphic P-diagrams for $W$. 
A Coxeter group $W$ is said to be {\it strongly rigid} 
if any two sets of Coxeter generators for $W$ are conjugate. 

A Coxeter system $(W,S)$ is said to be {\it complete} if 
the underlying graph of the P-diagram of $(W,S)$ is complete. 
A Coxeter system $(W,S)$ is said to be {\it finite} (resp.  
{\it infinite}) if $W$ is finite (resp. infinite).

\begin{theorem}{\rm (Caprace, Franzsen, 
Howlett, and M\"uhlherr \cite{C-M},\cite{F-H-M})} 
If $(W,S)$ is an infinite, complete, irreducible Coxeter system of finite rank, 
then $W$ is strongly rigid. 
\end{theorem}

\section{Coxeter Systems of Finite Coxeter Groups}  

We shall use Coxeter's notation on p. 297 of \cite{Coxeterb}
for the irreducible spherical Coxeter simplex reflection groups except 
that we denote the dihedral group ${\bf D}_2^k$ by ${\bf D}_2(k)$. 
Subscripts denote the rank of a Coxeter system in Coxeter's notation. 
Coxeter's notation partly agrees with but differs from Bourbaki's notation on p.193 of 
\cite{Bourbaki}.

Coxeter \cite{Coxeter} proved that every finite irreducible Coxeter system 
is isomorphic to exactly one 
of the Coxeter systems ${\bf A}_n$, ${\bf B}_n$, ${\bf C}_n$, ${\bf D}_2(k)$, 
${\bf E}_6$, ${\bf E}_7$, ${\bf E}_8$, ${\bf F}_4$, ${\bf G}_3$, ${\bf G}_4$ described below. 
Each of these Coxeter groups, of rank $n$, is 
a finite group of orthogonal $n\times n$ matrices. 
The center of each of these Coxeter groups is either $\{I\}$ or $\{\pm I\}$. 
We denote the center of a group $G$ by $Z(G)$. 
If $G$ is a group of orthogonal matrices, we denote 
the subgroup of determinant 1 matrices in $G$ by $G^+$. 

The {\it type} of a finite irreducible Coxeter system $(W,S)$ is the 
isomorphism type of $(W,S)$ represented by one of the systems 
${\bf A}_n$, ${\bf B}_n$, ${\bf C}_n$, ${\bf D}_2(k)$, 
${\bf E}_6$, ${\bf E}_7$, ${\bf E}_8$, ${\bf F}_4$, ${\bf G}_3$, ${\bf G}_4$. 
The {\it type} of an irreducible subset $A$ of $S$  
is the type of $(\langle A\rangle,A)$.  

The Coxeter group ${\bf A}_n$ is the group of symmetries of a regular $n$-simplex 
for each $n\geq 1$,  
and so ${\bf A}_n$ is isomorphic to the symmetric group $S_{n+1}$ for each $n\geq 1$. 
The C-diagram of ${\bf A}_n$ is the following linear diagram 
with $n$ vertices and all edge labels 3:

$$\mbox{
\setlength{\unitlength}{1cm}
\begin{picture}(4,.5)(0,-.3)
\thicklines
\put(0,0){\line(1,0){1}}
\put(1,0){\line(1,0){1}}
\put(2.3,0){$\ldots$}
\put(3,0){\line(1,0){1}}
\put(0,0){\circle*{.15}}
\put(1,0){\circle*{.15}}
\put(2,0){\circle*{.15}}
\put(3,0){\circle*{.15}}
\put(4,0){\circle*{.15}}
\end{picture}}$$
The Coxeter generators $a_1,\ldots,a_n$ of ${\bf A}_n$, indexed so that $m(a_i,a_{i+1})=3$
for $i=1,\ldots,n$, correspond to the transpositions $(12), (23),\ldots, (n\,n+1)$ of $S_{n+1}$.
The group ${\bf A}_n$ has order $(n+1)!$ for all $n\geq 1$. 
The center of ${\bf A}_n$ is trivial for all $n\geq 2$. 

The Coxeter group ${\bf C}_n$ is the group of symmetries 
of an $n$-cube for each $n\geq 2$,  
and ${\bf C}_n$ is represented by the group of all $n\times n$ 
orthogonal matrices in which each column has all zero entries except for one, 
which is $\pm 1$. Thus we have a split short exact sequence
$$1\lra D_n \lra {\bf C}_n \ {\buildrel \pi \over \lra }\ S_n \lra 1$$
where $D_n=\{\diag(\pm 1,\pm 1,\ldots,\pm 1)\}$ and $\pi$ maps 
a permutation matrix to the corresponding permutation. 
The C-diagram of ${\bf C}_n$ is the following linear diagram with $n$ vertices:

$$\mbox{
\setlength{\unitlength}{1cm}
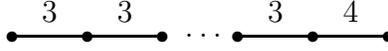
\begin{picture}(5,.8)(0,-.3)
\thicklines
\put(0,0){\line(1,0){1}}
\put(1,0){\line(1,0){1}}
\put(2.3,0){$\ldots$}
\put(3,0){\line(1,0){1}}
\put(4,0){\line(1,0){1}}
\put(0,0){\circle*{.15}}
\put(1,0){\circle*{.15}}
\put(2,0){\circle*{.15}}
\put(3,0){\circle*{.15}}
\put(4,0){\circle*{.15}}
\put(5,0){\circle*{.15}}
\put(.4,.2){3}
\put(1.4,.2){3}
\put(3.4,.2){3}
\put(4.4,.2){4}
\end{picture}}$$
The Coxeter generators $c_1,\ldots,c_n$ of ${\bf C}_n$ are indexed so that 
$m(c_i,c_{i+1})=3$ for $i=1,\ldots, n-2$ and $m(c_{n-1},c_n)=4$. 
The generators $c_1,\ldots,c_{n-1}$ are represented by the permutation 
matrices corresponding to the transpositions $(12), (23),\ldots, (n-1\,n)$ 
and $c_n$ is represented by the matrix $\diag(1,\ldots,1,-1)$. 
The order of the group ${\bf C}_n$ is $2^nn!$ and $Z({\bf C}_n) = \{\pm I\}$.

The Coxeter group ${\bf B}_n$, with $n\geq 4$, is a subgroup of ${\bf C}_n$ of index 2  
with Coxeter generators $b_i=c_i$, for $i=1,\ldots, n-1$, and $b_n = c_nc_{n-1}c_n$. 
We have $b_{n-1}b_n = \diag(1,\ldots,1,-1,-1)$ and  $m(b_{n-2},b_n) = 3$. 
The group ${\bf B}_n$ contains $D_n^+$ 
and the group of permutation matrices, and so we have a split 
short exact sequence
$$1\lra D_n^+ \lra {\bf B}_n \ {\buildrel \pi \over \lra }\ S_n \lra 1.$$
The C-diagram of ${\bf B}_n$ is the following Y-shaped diagram 
with $n$ vertices and all edge labels 3:

$$\mbox{
\setlength{\unitlength}{1cm}
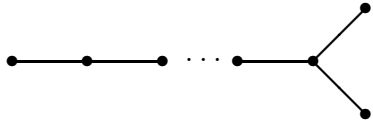
\begin{picture}(5,1.2)(0,-.6)
\thicklines
\put(0,0){\line(1,0){1}}
\put(1,0){\line(1,0){1}}
\put(2.3,0){$\ldots$}
\put(3,0){\line(1,0){1}}
\put(4,0){\line(1,1){.7}}
\put(4,0){\line(1,-1){.7}}
\put(0,0){\circle*{.15}}
\put(1,0){\circle*{.15}}
\put(2,0){\circle*{.15}}
\put(3,0){\circle*{.15}}
\put(4,0){\circle*{.15}}
\put(4.7,.7){\circle*{.15}}
\put(4.7,-.7){\circle*{.15}}
\end{picture}}$$

In order to have uniformity of notation, 
we extend the above definition of ${\bf B}_n$ to include the rank $n=3$. 
The group ${\bf B}_3$ is of type ${\bf A}_3$ and represents the degenerate 
case when there are no horizontal edges in the above diagram.   
The order of the group ${\bf B}_n$ is $2^{n-1}n!$ for each $n\geq 3$. 
The center of the group ${\bf B}_n$ is trivial if $n$ is odd 
and is $\{\pm I\}$ if $n$ is even. 

If $n \neq 4$, we call the two right most vertices, $b_{n-1}$ and $b_n$, 
of the above C-diagram of ${\bf B}_n$ the {\it split ends} of the diagram. 
We call any two endpoints of the C-diagram of ${\bf B}_4$ a pair 
of {\it split ends} of the diagram.

The group ${\bf D}_2(k)$ is the group of symmetries 
of a regular $k$-gon for each $k\geq 5$. 
In order to have uniformity of notation, 
we extend the definition of ${\bf D}_2(k)$ to include the cases $k=3,4$, 
and so ${\bf D}_2(k)$ is a dihedral group of order $2k$ for each $k\geq 3$. 
Note that ${\bf D}_2(3)$ is of type ${\bf A}_2$ and ${\bf D}_2(4)$ is of type ${\bf C}_2$. 

The C-diagram of ${\bf D}_2(k)$ is an edge with label $k$:  
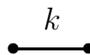
$$\mbox{
\setlength{\unitlength}{1cm}
\begin{picture}(1,.5)(0,0)
\thicklines
\put(0,0){\line(1,0){1}}
\put(0,0){\circle*{.15}}
\put(1,0){\circle*{.15}}
\put(.4,.25){$k$}
\end{picture}}$$
Let $a$ and $b$ be Coxeter generators for ${\bf D}_2(k)$. 
The center of ${\bf D}_2(k)$ is trivial if $k$ is odd 
and is generated by $(ab)^{k/2}$ if $k$ is even. 

The orders of ${\bf E}_6, {\bf E}_7, {\bf E}_8$ are 
$72\cdot 6!, 8\cdot 9!, 192\cdot 10!$, respectively. 
The center of ${\bf E}_6$ is trivial while the centers of ${\bf E}_7$ and ${\bf E}_8$ 
are $\{\pm I\}$. 
The C-diagrams of ${\bf E}_6$, ${\bf E}_7, {\bf E}_8$ are the following diagrams 
with all edge labels 3:

$$\mbox{
\setlength{\unitlength}{1cm}
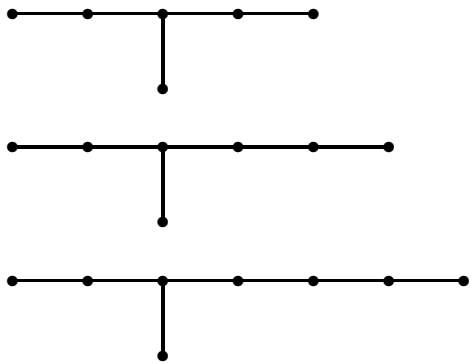
\begin{picture}(4,1.4)(0,-1)
\thicklines
\put(0,0){\line(1,0){1}}
\put(1,0){\line(1,0){1}}
\put(2,0){\line(1,0){1}}
\put(2,0){\line(0,-1){1}}
\put(3,0){\line(1,0){1}}
\put(0,0){\circle*{.15}}
\put(1,0){\circle*{.15}}
\put(2,0){\circle*{.15}}
\put(2,-1){\circle*{.15}}
\put(3,0){\circle*{.15}}
\put(4,0){\circle*{.15}}
\end{picture}}
$$

$$\mbox{
\setlength{\unitlength}{1cm}
\begin{picture}(4,1)(0,-1)
\thicklines
\put(0,0){\line(1,0){1}}
\put(1,0){\line(1,0){1}}
\put(2,0){\line(1,0){1}}
\put(2,0){\line(0,-1){1}}
\put(3,0){\line(1,0){1}}
\put(4,0){\line(1,0){1}}
\put(0,0){\circle*{.15}}
\put(1,0){\circle*{.15}}
\put(2,0){\circle*{.15}}
\put(2,-1){\circle*{.15}}
\put(3,0){\circle*{.15}}
\put(4,0){\circle*{.15}}
\put(5,0){\circle*{.15}}
\end{picture}}
$$

$$\mbox{
\setlength{\unitlength}{1cm}
\begin{picture}(4,1)(0,-1)
\thicklines
\put(0,0){\line(1,0){1}}
\put(1,0){\line(1,0){1}}
\put(2,0){\line(1,0){1}}
\put(2,0){\line(0,-1){1}}
\put(3,0){\line(1,0){1}}
\put(4,0){\line(1,0){1}}
\put(5,0){\line(1,0){1}}
\put(0,0){\circle*{.15}}
\put(1,0){\circle*{.15}}
\put(2,0){\circle*{.15}}
\put(2,-1){\circle*{.15}}
\put(3,0){\circle*{.15}}
\put(4,0){\circle*{.15}}
\put(5,0){\circle*{.15}}
\put(6,0){\circle*{.15}}
\end{picture}}
$$

The Coxeter group ${\bf F}_4$ is the group of symmetries of a regular 24-cell. 
The order of ${\bf F}_4$ is 1152 and $Z({\bf F}_4) = \{\pm I\}$. 
The C-diagram of ${\bf F}_4$ is the linear diagram: 

$$\mbox{
\setlength{\unitlength}{1cm}
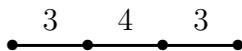
\begin{picture}(4,.8)(0,-.3)
\thicklines
\put(0,0){\line(1,0){1}}
\put(1,0){\line(1,0){1}}
\put(2,0){\line(1,0){1}}
\put(0,0){\circle*{.15}}
\put(1,0){\circle*{.15}}
\put(2,0){\circle*{.15}}
\put(3,0){\circle*{.15}}
\put(.4,.2){3}
\put(1.4,.2){4}
\put(2.4,.2){3}
\end{picture}}$$

The Coxeter group ${\bf G}_3$ is the group of symmetries of a regular dodecahedron. 
The order of ${\bf G}_3$ is 120 and the center of ${\bf G}_3$ has order two. 
The Coxeter group ${\bf G}_4$ is the group of symmetries of a regular 120-cell.  
The order of ${\bf G}_4$ is $120^2$ and $Z({\bf G}_4)= \{\pm I\}$. 
The C-diagrams of ${\bf G}_3$ and ${\bf G}_4$ are the linear diagrams:

$$\mbox{
\setlength{\unitlength}{1cm}
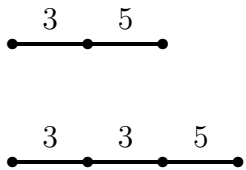
\begin{picture}(4,.8)(0,-.3)
\thicklines
\put(0,0){\line(1,0){1}}
\put(1,0){\line(1,0){1}}
\put(0,0){\circle*{.15}}
\put(1,0){\circle*{.15}}
\put(2,0){\circle*{.15}}
\put(.4,.2){3}
\put(1.4,.2){5}
\end{picture}}$$

$$\mbox{
\setlength{\unitlength}{1cm}
\begin{picture}(4,.8)(0,-.3)
\thicklines
\put(0,0){\line(1,0){1}}
\put(1,0){\line(1,0){1}}
\put(2,0){\line(1,0){1}}
\put(0,0){\circle*{.15}}
\put(1,0){\circle*{.15}}
\put(2,0){\circle*{.15}}
\put(3,0){\circle*{.15}}
\put(.4,.2){3}
\put(1.4,.2){3}
\put(2.4,.2){5}
\end{picture}}$$

Lemmas 3.1 through 3.7 are either elementary or well known. 

\begin{lemma}  
The Coxeter groups ${\bf A}_n$ and ${\bf B}_n$ are indecomposable with respect 
to direct products for all $n$. 
\end{lemma}

\begin{lemma} 
The Coxeter group ${\bf C}_n$ 
is decomposable with respect to direct products 
if and only if $n$ is odd.  
If $n$ is odd and ${\bf C}_n= H\times K$ with $1 < |H|\leq |K|$, 
then $H=\{\pm I\}$ and $K = {\bf B}_n$ or $\theta({\bf B}_n)$ 
where $\theta$ is the automorphism of ${\bf C}_n$ defined 
by $\theta(c_i)=-c_i$, for $i=1,\ldots,n-1$, and $\theta(c_n)=c_n$. 
\end{lemma} 

\begin{lemma}  
The Coxeter group ${\bf D}_2(n)$, with Coxeter generators $a$ and $b$, 
is decomposable with respect to direct products 
if and only if $n \equiv 2\ {\rm mod}\ 4$.  
If $n \equiv 2\ {\rm mod}\ 4$ and ${\bf D}_2(n) = H\times K$ with $1 < |H|\leq |K|$, 
then $H=\langle (ab)^{n/2}\rangle$ and $K = \langle a, bab \rangle$ or $\langle b, aba\rangle$, 
moreover $K \cong {\bf D}_2(n/2)$.  
\end{lemma} 

\begin{lemma}  
The Coxeter groups ${\bf E}_6$ and ${\bf E}_8$ are indecomposable with respect to direct products. 
\end{lemma}

\begin{lemma}  
The Coxeter group ${\bf E}_7$ is decomposable with respect to direct products. 
If ${\bf E}_7 = H \times K$ with $1 < |H|\leq |K|$, 
then $H = \{\pm I\}$ and $K = {\bf E}_7^+$, 
moreover ${\bf E}_7^+$ is a nonabelian simple group. 
\end{lemma}

\begin{lemma} 
The Coxeter groups ${\bf F}_4$  and ${\bf G}_4$ are indecomposable with respect to direct products. 
\end{lemma}

\begin{lemma} 
The Coxeter group ${\bf G}_3$ is decomposable with respect to direct products. 
If ${\bf G}_3 = H \times K$ with $1 < |H|\leq |K|$, 
then $H = \{\pm I\}$ and $K = {\bf G}_3^+$, moreover ${\bf G}_3^+$ is a nonabelian simple group. 
\end{lemma}

The next lemma follows from the Krull-Remak-Schmidt Theorem (KRS-Theorem), 
Theorem 4.8 in \cite{Suzuki}.
\begin{lemma} 
Let $G$ be a finite group with direct product decompositions
$$G = H_1 \times H_2 \times \cdots \times H_r
\quad\hbox{and}\quad
G = K_1 \times K_2 \times \cdots \times K_s$$
such that $H_i$ and $K_j$ are nontrivial and indecomposable 
with respect to direct products for each $i$ and $j$. 
Let $\iota_i:H_i\to G$ be the inclusion map for each $i$ 
and let $\pi_j:G\to K_j$ be the projection map for each $j$. 
Suppose $H_p$ is nonabelian. 
Then there is a unique $q$ such that 
$H_p\cap K_q \not=\{1\}$.  Moreover $\pi_q\iota_p:H_p\to K_q$ 
is an isomorphism and $Z(G)H_p = Z(G)K_q$. 
Furthermore $[H_p,H_p]=[K_q,K_q]$ and $\pi_q\iota_p:H_p\to K_q$ 
restricts to the identity on $[H_p,H_p]$. 
\end{lemma}

\begin{theorem} {\rm (Matching Theorem for Systems of a Finite Coxeter Group)} 
Let $W$ be a finite Coxeter group with two sets of Coxeter generators $S$ and $S'$. 
Let 
$$(W,S) = (W_1,S_1)\times \cdots\times (W_m,S_m)$$
and
$$(W,S') = (W_1',S_1')\times \cdots\times (W_n',S_n')$$
be the factorizations of $(W,S)$ and $(W,S')$ into irreducible factors. 
Let $k$ be such that $W_k$ is noncyclic. 
Then there is a unique $\ell$ such that $W_\ell'$ is noncyclic and 
$[W_k,W_k]=[W_\ell',W_\ell']$. 
Moreover, 
\begin{enumerate}
\item $Z(W)W_k = Z(W)W_\ell'$.
\item If $|W_k| = |W_\ell'|$, then $(W_k,S_k)\cong (W_\ell',S_\ell')$ and 
there is an isomorphism $\phi: W_k\to W_\ell'$ that restricts 
to the identity on $[W_k,W_k]$. 
\item If $|W_k| < |W_\ell'|$, then either $(W_k,S_k)$ has type ${\bf B}_{2q+1}$ and 
$(W_\ell',S_\ell')$ has type ${\bf C}_{2q+1}$ for some $q\geq 1$ or 
$(W_k,S_k)$ has type ${\bf D}_2(2q+1)$ and 
$(W_\ell',S_\ell')$ has type ${\bf D}_2(4q+2)$ for some $q\geq 1$,  
and there is a monomorphism $\phi: W_k \to W_\ell'$ that restricts 
to the identity on $[W_k,W_k]$.
\end{enumerate}
\end{theorem}
\begin{proof}
By Lemmas 3.1-3.7, we can refine the decomposition $W=W_1\times \cdots \times W_m$ 
to a decomposition $W=H_1\times\cdots\times H_r$, 
with $H_i$ nontrivial and indecomposable 
with respect to direct products, by replacing each $W_i$ that factors into 
a direct product $W_i = H_{j-1}\times H_j$, with $|H_{j-1}|=2$, by $H_{j-1}\times H_j$. 
Likewise refine the decomposition $W=W_1'\times \cdots \times W_\ell'$ to a decomposition 
$W=K_1\times\cdots\times K_s$, with $K_i$ nontrivial and indecomposable 
with respect to direct products, by replacing each $W_i'$ that factors into 
a direct product $W_i' = K_{j-1}\times K_j$, with $|K_{j-1}|=2$, by $K_{j-1}\times K_j$. 
Then $r=s$ by the KRS-Theorem. 

Suppose that $W_k$ is noncyclic. Then $W_k$ is nonabelian, 
since $(W_k,S_k)$ is irreducible. Now $W_k = H_p$ 
or $H_{p-1}\times H_p$, with $|H_{p-1}|=2$, for some $p$. 
In either case $H_p$ is nonabelian by Lemmas 3.1-3.7. 
By Lemma 3.8, there is a unique $q$ such that $H_p\cap K_q\neq\{1\}$. 
Moreover $[H_p,H_p]=[K_q,K_q]$ and $\pi_q\iota_p:H_p\to K_q$ 
restricts to the identity on $[H_p,H_p]$. Then $K_q$ is nonabelian. 
Now there is an $\ell$ such that $W_\ell'=K_q$ or $K_{q-1}\times K_q$ with $|K_{q-1}|=2$. 
Then $W_\ell'$ is noncyclic and 
$$[W_k,W_k]= [H_p,H_p]=[K_q,K_q]=[W_\ell',W_\ell'].$$

Now suppose $W_i'$ is noncyclic and $[W_k,W_k]=[W_i',W_i']$. 
Then $W_i'$ is nonabelian, since $W_k$ is nonabelian. 
Now $W_i' = K_j$ or $K_{j-1}\times K_j$, with $|K_{j-1}| = 2$, for some $j$. 
Then we have 
$$[H_p,H_p]= [W_k,W_k]=[W_i',W_i']=[K_j,K_j].$$
Hence $H_p\cap K_j\neq\{1\}$, and so $j=q$ by the uniqueness of $q$. 
Therefore $K_q\subset W_\ell'\cap W_i'$, 
and so $i=\ell$ and $\ell$ is unique. 

(1)  By Lemma 3.8, we have
$\!Z(W)W_k = Z(W)H_p = Z(W)H_q = Z(W)W'_\ell.$

(2)  Suppose $|W_k|=|W_\ell'|$.  As $H_p\cong H_q$, we have either 
$W_k = H_p$ and $W_\ell' =K_q$ or $W_k=H_{p-1}\times H_p$, with $|H_{p-1}|=2$, 
and $W_\ell'=K_{q-1}\times K_q$, with $|K_{q-1}|=2$. Hence $W_k\cong W_\ell'$ 
and $(W_k,S_k)\cong (W_\ell',S_\ell')$, since $(W_k,S_k)$ and $(W_\ell',S_\ell')$ 
are irreducible. 
Moreover $\pi_q\iota_p:H_p\to K_q$ is an isomorphism that 
restricts to the identity on $[H_p,H_p]$. 
If $W_k=H_p$ and $W_\ell'=K_q$, let $\phi=\pi_q\iota_p$. 
If $W_k=H_{p-1}\times H_p$ and $W_\ell'=H_{q-1}\times H_q$ extend 
$\pi_q\iota_p:H_p\to K_q$ to an isomorphism $\phi:W_k\to W_\ell'$ by mapping 
the generator of $H_{p-1}$ to the generator of $H_{q-1}$. 
Then $\phi:W_k\to W_\ell'$ is an isomorphism that restricts to the identity 
on $[W_k,W_k]=[H_p,H_p]$.

(3) Suppose $|W_k|<|W_\ell'|$. As $H_p\cong H_q$, we have 
$W_k = H_p$ and $W_\ell'=K_{q-1}\times K_q$, with $|K_{q-1}|=2$.
By Lemmas 3.1-3.7, either $(W_k,S_k)\cong {\bf B}_{2q+1}$ and 
$(W_\ell',S_\ell') \cong {\bf C}_{2q+1}$ for some $q\geq 1$ or 
$(W_k,S_k)\cong {\bf D}_2(2q+1)$ and 
$(W_\ell',S_\ell') \cong {\bf D}_2(4q+2)$ for some $q\geq 1$. 
Moreover $\pi_q\iota_p:H_p\to K_q$ is an isomorphism that 
restricts to the identity on $[H_p,H_p]$. 
Hence $\pi_q\iota_p:H_p\to K_q$ extends to a monomorphism $\phi:W_k\to W_\ell'$ 
that restricts to the identity on $[W_k,W_k]$.
\end{proof}

\section{The Basic Matching Theorem}

Let $(W,S)$ be a Coxeter system. 
The undirected {\it Cayley graph} of $(W,S)$ is graph 
${\rm K} = {\rm K}(W,S)$ with vertices $W$ 
and edges unordered pairs $(v,w)$ such that $w=vs$ for some element $s$ of $S$. 
The element $s=v^{-1}w$ of $S$ is called the {\it label} of the edge $(v,w)$. 
We represent an edge path in ${\rm K}$ 
beginning at vertex $v$ by ``$\alpha = (s_1,\ldots ,s_n)$ at $v$" 
where $s_i$ is the label of the $i$th edge of the path. 
The {\it length of an edge path} $\alpha = (s_1,\ldots ,s_n)$ is $|\alpha| = n$. 
The {\it distance} between vertices $v$ and $w$ in ${\rm K}$ is the minimal length $d(v,w)$ 
of an edge path from $v$ to $w$. 
A {\it geodesic} in ${\rm K}$ is an edge path $\alpha$ 
from a vertex $v$ to a vertex $w$ such that $|\alpha|= d(v,w)$ in ${\rm K}$.  
The {\it length of an element} $w$ of $W$ is $l(w) = d(1,w)$. 
A word $w=s_1s_2\cdots s_n$, with $s_i$ in $S$, is said to be {\it reduced} if $n=l(w)$. 

\begin{lemma} 
If $A\subset S$, then for any pair of vertices $v , w$ in ${\rm K}$, 
there is a unique element $x$ of the coset $w\langle A\rangle$
nearest to $v$ and for any geodesic $\alpha$ from $v$ to $x$ 
and geodesic $\beta$ at $x$ in $w\langle A\rangle$ 
(i.e. all edge labels of $\beta$ are in $A$), 
the path $\alpha \beta$ is geodesic. 
Moreover an element $x$ of $w\langle A\rangle$ is the nearest element 
of $w\langle A\rangle$ to $v$ if and only if for any geodesic $\alpha$ from $v$ to $x$ 
the path $(\alpha,a)$ is a geodesic for each $a$ in $A$. 
\end{lemma}
\begin{proof} 
Suppose $x$ and $y$ are distinct elements of $w\langle A\rangle$ that are nearest to $v$. 
Let $\alpha$ and $\gamma$ be geodesics from $v$ to $x$ and $y$, respectively. 
Then $|\alpha| = |\gamma|$. 
Let $\beta$ be a geodesic, with labels in $A$, from $x$ to $y$. 
The path $\alpha\beta$ is not geodesic, since $|\alpha\beta| > |\gamma|$,  
and so a letter of $\beta$ deletes with a letter of $\alpha$ by the deletion condition. 
This defines a path from $v$ to $w\langle A\rangle$ shorter than $\alpha$, 
which is impossible. 
A proof of the second assertion of the lemma is analogous. 

Now suppose $y$ is an element of $w\langle A\rangle$ such that 
for any geodesic $\gamma$ from $v$ to $y$ the path $(\gamma,a)$ is geodesic for each $a$ in $A$. 
Then $y$ is the nearest element $x$ of $w\langle A\rangle$ to $v$ 
otherwise there would be a geodesic $\alpha\beta$ from $v$ to $y$ 
with $\alpha$ a geodesic from $v$ to $x$ and $\beta$ a nontrivial 
geodesic, with labels in $A$, from $x$ to $y$,  
but $\beta$ ends in some element $a$ of $A$, and so 
the path $(\alpha\beta,a)$ would not be geodesic. 
\end{proof}

\begin{lemma} {\rm (Bourbaki \cite{Bourbaki}, Ch. IV, \S 1, Ex. 3)} 
If $A,B\subset S$ and $w$ is an element of $W$, then there is a unique shortest
representative $x$ of the double coset $\langle A\rangle w\langle B\rangle$.
\end{lemma}

\begin{lemma} 
Let $A,B\subset S$ and let $w$ in $W$ be such that 
$w\langle A\rangle w^{-1}\subset \langle B\rangle$.  
If $u$ is the shortest element of 
$\langle B\rangle w\langle A\rangle$, then $uAu^{-1}\subset B$.
\end{lemma}
\begin{proof} Certainly we have $u\langle A\rangle u^{-1}\subset\langle B\rangle$. 
Let $u = u_1\cdots u_n$ be reduced.
For any $a$ in $A$, the word $ua= u_1\cdots u_na$ is reduced by Lemma 4.1. 
Now $uau^{-1}$ is in $\langle B\rangle$.  Write $uau^{-1} = b_1\cdots b_k$ 
with $b_1\cdots b_k$ reduced in $\langle B\rangle$. 
Now $u$ is a shortest element of $\langle B\rangle u$, 
and so $u^{-1}$ is the shortest element of $u^{-1}\langle B\rangle$. 
Hence $u_n\cdots u_1b_k\cdots b_1$ is reduced by Lemma 4.1, 
and so $b_1\cdots b_ku_1\cdots u_n$ is reduced. 
As $u_1\cdots u_na=b_1\cdots b_ku_1\cdots u_n$, 
we have $k=1$ and $uau^{-1}=b_1$. 
\end{proof}

\begin{lemma} {\rm (Bourbaki \cite{Bourbaki}, Ch. IV, \S 1, Ex. 22)} 
Let $w_0$ be an element of $W$. Then the following are equivalent.
\begin{enumerate}
\item $l(w_0s)<l(w_0)$ for all $s$ in $S$.
\item $l(w_0w)=l(w_0)-l(w)$ for all $w$ in $W$.
\end{enumerate}
Such an element $w_0$ is unique and exists if and only if $W$ is finite. 
If $W$ is finite, then $w_0$ is the unique element of maximal length in $W$. 
Moreover $w_0^2=1$ and $w_0Sw_0 = S$. 
\end{lemma} 

Let $(W,S)$ be a Coxeter system.  The {\it quasi-center} of $(W,S)$ is the subgroup 
$$QZ(W,S) =\{w \in W: wSw^{-1}=S\}.$$

\begin{lemma} {\rm (Bourbaki \cite{Bourbaki}, Ch. V, \S 4, Ex. 3)} 
Let $(W,S)$ be an irreducible Coxeter system with a nontrivial quasi-center. 
Then $W$ is a finite group and $QZ(W)=\{1,w_0\}$ with $w_0$ the longest element of $(W,S)$. 
\end{lemma} 

Let $V$ be a real vector space having a basis $\{e_s: s\in S\}$ 
in one-to-one correspondence with $S$. 
Let $B$ be the symmetric bilinear form on $V$ defined by
$$B(e_s,e_t) =
\begin{cases}
 -\cos(\pi/m(s,t)) & {\rm if}\ m(s,t)<\infty,\\
               -1 & {\rm if}\ m(s,t)=\infty.
\end{cases}$$
There is an action of $W$ on $V$ defined by 
$$s(x)=x-2B(x,e_s)e_s \quad\hbox{for all}\ s\in S\ \hbox{and}\ x\in V.$$

The {\it root system} of $(W,S)$ is the set 
of vectors 
$$\Phi = \{w(e_s): w\in W\ \hbox{and}\ s\in S\}.$$
The elements of $\Phi$ are called {\it roots}. 
By Prop. 2.1 of Deodhar \cite{Deodhar}, every root $\phi$ can be written uniquely in the form 
$\phi = \sum_{s\in S}a_se_s$ with $a_s\in \realnos$ 
where either $a_s\geq 0$ for all $s$ or $a_s\leq 0$ for all $s$. 
In the former case, we say $\phi$ is {\it positive} and write $\phi>0$.  
Let $\Phi^+$ be the set of positive roots. 

The set of {\it reflections} of $(W,S)$ is the set 
$$T=\{wsw^{-1}: w\in W\ \hbox{and}\ s\in S\}.$$
\begin{proposition} {\rm (Deodhar \cite{Deodhar}, Prop. 3.1)} 
The function $\rho: \Phi^+ \to T$ defined by $\rho(w(e_s))=wsw^{-1}$ 
is well defined and a bijection.
\end{proposition}
\begin{proposition} {\rm (Deodhar \cite{Deodhar}, Prop. 2.2)} 
Let $w \in W$ and $s\in S$. 
Then $l(ws) > l(w)$ if and only if $w(e_s) > 0$. 
\end{proposition}

If $A \subset S$, set $E_A =\{e_s: s\in A\}$. 
The next lemma follows from Lemma 4.1 and Propositions 4.6 and 4.7. 

\begin{lemma} 
Let $A, B \subset S$ and let $w\in W$. 
Then the following are equivalent: 
\begin{enumerate}
\item $w(E_A)=E_B$. 
\item $wAw^{-1}=B$ and $l(wa) > l(w)$ for all $a \in A$. 
\item $wAw^{-1}=B$ and $w$ is the shortest element of $w\langle A\rangle$. 
\end{enumerate}
\end{lemma}

The next lemma follows from Lemma 4.3 and Lemma 4.8. 
\begin{lemma} 
Let $A,B\subset S$ and let $w$ in $W$ be such that 
$w\langle A\rangle w^{-1} = \langle B\rangle$.  
If $u$ is the shortest element of 
$\langle B\rangle w\langle A\rangle$, then $u(E_A)= E_B$.
\end{lemma}

Suppose $A\subset S$. If $\langle A\rangle$ is finite, 
we denote the longest element of $\langle A\rangle$ by $w_A$. 
Suppose $s\in S-A$.  
Let $K\subset S$ be the irreducible component of $A\cup\{s\}$ containing $s$. 
We say that $s$ is $A$-{\it admissible} if $\langle K\rangle$ is finite. 
If $s$ is $A$-admissible, define 
$$\nu(s,A)=w_Kw_{K-\{s\}}.$$
Then $\nu(s,A)$ is the shortest element of $w_K\langle A\rangle$ by Lemma 4.4;   
moreover, if $t = w_Ksw_K$ and $B = (A\cup \{s\})-\{t\}$,  
then $\nu(s,A)(E_A)= E_B$ by Lemma 4.8. 

\begin{proposition} {\rm (Deodhar \cite{Deodhar}, Prop. 5.5)} 
Let $A, B\subset S$, and let $w\in W$.  If $w(E_A) = E_B$ and $w\neq 1$, 
then there exists a sequence $A_1,A_2,\ldots, A_{n+1}$ of subsets of $A$,
and a sequence $s_1,s_2,\ldots, s_n$ 
of elements of $S$ such that 
\begin{enumerate}
\item $A_1=A$ and $A_{n+1}=B$,
\item $s_i \in S-A_i$ and $s_i$ is $A_i$-admissible for $i=1,\ldots,n$,
\item $\nu(s_i,A_i)(E_{A_i})=E_{A_{i+1}}$ for $i=1,\ldots,n$, 
\item $w=\nu(s_n,A_n)\cdots \nu(s_2,A_2)\nu(s_1,A_1)$,
\item $l(w)=l(\nu(s_1,A_1))+l(\nu(s_2,A_2))+\cdots + l(\nu(s_n,A_n))$.
\end{enumerate}
\end{proposition}

The next lemma follows from Proposition 4.10.

\begin{lemma} 
Let $A\subset S$. Then there exists $B\subset S$ such that $A\neq B$ 
and $\langle A\rangle$ is conjugate to $\langle B\rangle$ in $W$ if and only if 
there exists $s\in S-A$ such that 
\begin{enumerate}
\item $m(s,a)> 2$ for some $a\in A$, 
\item the element $s$ is $A$-admissible,
\item if $K$ is the component of $A\cup\{s\}$ containing $s$, 
then $w_Ksw_K\neq s$. 
\end{enumerate}
\end{lemma}

\begin{lemma} 
Let $A, B \subset S$. 
If $\langle A\rangle$ is a maximal finite visual subgroup of $(W,S)$ 
and $\langle A\rangle$ and $\langle B\rangle$ are conjugate, then $A = B$. 
\end{lemma}
\begin{proof}
If $s\in S-A$, then the irreducible component of $\langle A\cup\{s\}\rangle$ 
containing $s$ is infinite, since $\langle A\rangle$ is a maximal finite visual 
subgroup of $(W,S)$.  
Hence no $s\in S-A$ is $A$-admissible, and so $A=B$ by Lemma 4.11.
\end{proof}

\begin{proposition} {\rm (Bourbaki \cite{Bourbaki}, Ch. V, \S 4. Ex. 2)} 
If $H$ is a finite subgroup of $W$, 
then there is a subset $A$ of $S$ such that $\langle A\rangle$ is finite 
and $H$ is conjugate to a subgroup of $\langle A\rangle$. 
\end{proposition}

\begin{lemma} 
Every maximal finite visual subgroup of $(W,S)$ is a maximal finite subgroup of $W$. 
\end{lemma}
\begin{proof}
Let $M\subset S$ be such that $\langle M\rangle$ 
is a maximal finite visual subgroup of $(W,S)$. 
Suppose $H$ is a finite subgroup of $W$ containing $\langle M\rangle$. 
Then $wHw^{-1}\subset \langle A\rangle$ for some $w\in W$ 
and some $A\subset S$ such that $\langle A\rangle$ is finite by Prop. 4.13. 
Then $w\langle M\rangle w^{-1}\subset \langle A\rangle$. 
Let $u$ be the shortest element of $\langle A\rangle w\langle M\rangle$. 
Then $uMu^{-1}\subset A$ by Lemma 4.3. 
As no element of $S-M$ is $M$-admissible, $uMu^{-1}=M$ by Prop. 4.10. 
Therefore $M=A$, since $M$ is a maximal finite visual subgroup. 
Hence $w\langle M\rangle w^{-1} =\langle A\rangle$, and so $\langle M\rangle = H$. 
Thus $\langle M\rangle$ is a maximal finite subgroup of $W$. 
\end{proof}

A {\it simplex} $C$ of $(W,S)$ is a subset $C$ of $S$ such that 
$(\langle C\rangle, C)$ is a complete Coxeter system. 
A simplex $C$ of $(W,S)$ is said to be {\it spherical} 
if $\langle C\rangle$ is finite.  
The next proposition follows from Proposition 4.13 and Lemmas 4.12 and 4.14. 

\begin{proposition}  
Let $W$ be a finitely generated Coxeter group with two sets 
of Coxeter generators $S$ and $S'$, 
and let $M$ be a maximal spherical simplex of $(W,S)$. 
Then there is a unique maximal spherical simplex $M'$ of $(W,S')$ 
such that $\langle M\rangle$ and $\langle M'\rangle$ 
are conjugate in $W$. 
\end{proposition} 

The next lemma follows from Lemma 4.11. 

\begin{lemma} 
Let $A, B \subset S$. 
If $\langle A\rangle$ is a maximal finite irreducible subgroup of $(W,S)$ 
and $\langle A\rangle$ and $\langle B\rangle$ are conjugate, then $A = B$. 
\end{lemma}

\begin{lemma} 
Let $x,y\in S$ be distinct, let $B\subset S$,  
and let $w\in W$ such that $wxyw^{-1}\in \langle B\rangle$. 
If $u$ is the shortest element of the
double coset $\langle B\rangle w\langle x,y\rangle$, then
$u\{x,y\}u^{-1}\subset B$. 
\end{lemma}
\begin{proof}  
Let $b = uxyu^{-1}$. Then $b$ is in $\langle B\rangle$. 
Write $b=b_1\cdots b_k$ with $b_1\cdots b_k$ reduced in $\langle B\rangle$ 
and let $u=u_1\cdots u_n$ be reduced.  
Observe that the left and right sides of  
$b_1\cdots b_ku_n\cdots u_1=u_1\cdots u_nxy$
are reduced, and so $k=2$. 
Consider the bigon in ${\rm K}$ with geodesic sides
$(u_1,\ldots , u_n,x,y)$ and $(b_1,b_2,u_1,\ldots ,u_n)$. 
The path $(b_1,u_1,\ldots ,u_n,x,y)$ is not geodesic. 
If $b_1$ deletes with $x$, then $uxu^{-1}=b_1$; 
otherwise, $b_1$ deletes with $y$ and $uxu^{-1}=b_2$. 
Now $(u_1,\ldots ,u_n,y,x)$ is geodesic, 
but $(b_2,u_1,\ldots ,u_n,y,x)$ is not. 
If $b_2$ deletes with $y$, then $uyu^{-1}=b_2$; 
otherwise, $b_2$ deletes with $x$ and $uyu^{-1}=b_1$. 
\end{proof}

\begin{lemma} 
Let $A,B\subset S$ with no $a\in A$ central in $\langle A\rangle$. 
Suppose 
$$w[\langle A\rangle, \langle A\rangle]w^{-1}\subset \langle B\rangle
\quad\hbox{for some}\ w\in W.$$
Let $u$ be the shortest element of $\langle B\rangle w \langle A\rangle$.  
Then $uAu^{-1}\subset B$.
\end{lemma}
\begin{proof}
Let $x\in A$. 
Then there exists $y\in A$ such that $m(x,y)>2$. 
Assume $m(x,y)$ is odd.  
Then $xy\in [\langle A\rangle,\langle A\rangle]$. 
Hence $uxu^{-1}\in B$ by Lemma 4.17. 

Assume now that $m(x,y)$ is even. 
Then $xyxy\in [\langle A\rangle,\langle A\rangle]$. 
Let $b= uxyxyu^{-1}$. Then $b$ is in $\langle B\rangle$. 
Write $b=b_1\cdots b_k$ with $b_1\cdots b_k$ reduced in $\langle B\rangle$ 
and let $u=u_1\cdots u_n$ be reduced.  
Observe that the left and right sides of 
$u_1\cdots u_nxyxy=b_1\cdots b_ku_1\cdots u_n$
are reduced, and so $k=4$. 

Assume that $m(x,y) > 4$. 
Consider the bigon in ${\rm K}$ with geodesic sides
$(u_1,\ldots,u_n,x,y,x,y)$ and $(b_1,\ldots, b_4,u_1,\ldots,u_n)$. 
See Figure 1. 
The word $yxyxu_n\cdots u_1b_1$ is not reduced  
and $b_1$ must delete with one of the first four letters. 
Also $xyxyxu_n\cdots u_1$ is reduced, so $xyxyxu_n\cdots u_1b_1$ 
has length $n+4$. But then, in this last word, $b_1$ cannot delete
with the second, third, or fourth letter. 
Hence in the word, $yxyxu_n\cdots u_1b_1$, the letter $b_1$ must delete with
the fourth letter. 
This shows that $xu^{-1}b_1=u^{-1}$, and so
$uxu^{-1}=b_1$, as desired.

\medskip

$$\mbox{
\setlength{\unitlength}{.8cm}
\begin{picture}(10,5)(0,0)
\thicklines
\put(0,0){\line(0,1){5}}
\put(0,1){\line(1,0){10}}
\put(10,1){\line(0,1){4}}
\put(0,5){\line(1,0){10}}
\put(0,0){\circle*{.15}}
\put(0,1){\circle*{.15}}
\put(0,2){\circle*{.15}}
\put(0,3){\circle*{.15}}
\put(0,4){\circle*{.15}}
\put(0,5){\circle*{.15}}
\put(1,1){\circle*{.15}}
\put(1,5){\circle*{.15}}
\put(9,1){\circle*{.15}}
\put(9,5){\circle*{.15}}
\put(10,1){\circle*{.15}}
\put(10,2){\circle*{.15}}
\put(10,3){\circle*{.15}}
\put(10,4){\circle*{.15}}
\put(10,5){\circle*{.15}}
\put(-.45,.35){x}
\put(-.45,1.35){y}
\put(-.45,2.35){x}
\put(-.45,3.35){y}
\put(-.45,4.35){x}
\put(.35,.55){$u_n$}
\put(.35,5.35){$u_n$}
\put(9.35,.55){$u_1$}
\put(9.35,5.35){$u_1$}
\put(10.3,1.35){$b_4$}
\put(10.3,2.35){$b_3$}
\put(10.3,3.35){$b_2$}
\put(10.3,4.35){$b_1$}
\end{picture}}$$
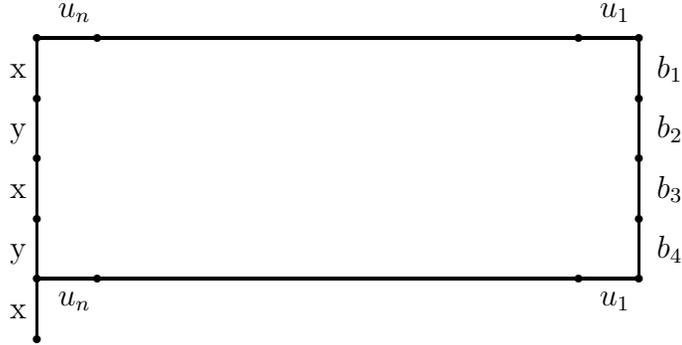
\centerline{\bf Figure 1}

\medskip

Assume now that $m(x,y)=4$. 
The word $yxyxu_n\cdots u_1b_1$ is not reduced  
and represents an element of length $n+3$. See Figure 1. 
In particular, $b_1$ does not delete with the
second or third letter of this word.  
Consequently, $b_1$ deletes with the first or fourth letter of the word. 
If $b_1$ deletes with the fourth letter, then
$uxu^{-1}=b_1$ and we are done. 
Suppose $b_1$ deletes with the first letter.   
Then $(uxyx)y(xyxu^{-1})=b_1$, and so $uyu^{-1}=b_1$. 
Now 
$$b_1b_2b_3b_4  =  uxyxyu^{-1} 
                =  uyxyxu^{-1} 
                =  uyu^{-1}uxyxu^{-1}  =  b_1uxyxu^{-1}.$$
%
Hence $uxyxu^{-1}=b_2b_3b_4$. 
Combining this last fact with the fact
that the word $u_1\cdots u_nxyxy$ is reduced allows us to use the
technique of the previous case to show that $uxu^{-1}=b_2$. 
See Figure 2.  
\end{proof}

\medskip

$$\mbox{
\setlength{\unitlength}{.8cm}
\begin{picture}(10,5)(0,0)
\thicklines
\put(0,0){\line(0,1){5}}
\put(0,1){\line(1,0){10}}
\put(0,4){\line(1,0){10}}
\put(10,1){\line(0,1){4}}
\put(0,5){\line(1,0){10}}
\put(0,0){\circle*{.15}}
\put(0,1){\circle*{.15}}
\put(0,2){\circle*{.15}}
\put(0,3){\circle*{.15}}
\put(0,4){\circle*{.15}}
\put(0,5){\circle*{.15}}
\put(1,1){\circle*{.15}}
\put(1,4){\circle*{.15}}
\put(1,5){\circle*{.15}}
\put(9,1){\circle*{.15}}
\put(9,4){\circle*{.15}}
\put(9,5){\circle*{.15}}
\put(10,1){\circle*{.15}}
\put(10,2){\circle*{.15}}
\put(10,3){\circle*{.15}}
\put(10,4){\circle*{.15}}
\put(10,5){\circle*{.15}}
\put(-.45,.35){y}
\put(-.45,1.35){x}
\put(-.45,2.35){y}
\put(-.45,3.35){x}
\put(-.45,4.35){y}
\put(.35,.55){$u_n$}
\put(.35,5.3){$u_n$}
\put(.35,4.3){$u_n$}
\put(9.35,.55){$u_1$}
\put(9.35,4.3){$u_1$}
\put(9.35,5.3){$u_1$}
\put(10.3,1.35){$b_4$}
\put(10.3,2.35){$b_3$}
\put(10.3,3.35){$b_2$}
\put(10.3,4.35){$b_1$}
\end{picture}}$$
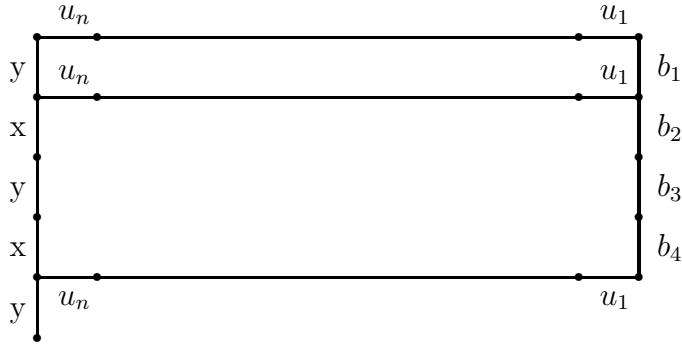
\centerline{\bf Figure 2}

\bigskip

Let $(W,S)$ be a Coxeter system. 
A {\it basic subgroup} of $(W,S)$ is a noncyclic, maximal, finite, irreducible, 
visual subgroup of $(W,S)$. 
A ${\it base}$ of $(W,S)$ is a subset $B$ of $S$ such that 
$\langle B\rangle$ is a basic subgroup of $(W,S)$.

\begin{theorem} {\rm (Basic Matching Theorem)} 
Let $W$ be a finitely generated Coxeter group with 
two sets of Coxeter generators $S$ and $S'$. 
Let $B$ be a base of $(W,S)$. 
Then there is a unique irreducible subset $B'$ of $S'$ such that 
$[\langle B\rangle,\langle B\rangle]$ is conjugate to 
$[\langle B'\rangle,\langle B'\rangle]$ in $W$. Moreover, 

\begin{enumerate}
\item The set $B'$ is a base of $(W,S')$, 

\item If $|\langle B\rangle|=|\langle B'\rangle|$, then $B$ and $B'$ have the same type 
and there is an isomorphism $\phi:\langle B\rangle \to \langle B'\rangle$ 
that restricts to conjugation on $[\langle B\rangle,\langle B\rangle]$ 
by an element of $W$. 

\item If $|\langle B\rangle|<|\langle B'\rangle|$, then either
$B$ has type ${\bf B}_{2q+1}$ and 
$B'$ has type ${\bf C}_{2q+1}$ for some $q\geq 1$ or 
$B$ has type ${\bf D}_2(2q+1)$ and 
$B'$ has type ${\bf D}_2(4q+2)$ for some $q\geq 1$. 
Moreover, there is a monomorphism $\phi:\langle B\rangle \to \langle B'\rangle$ 
that restricts to conjugation on $[\langle B\rangle,\langle B\rangle]$ 
by an element of $W$. 
\end{enumerate}
\end{theorem}
\begin{proof}
Let $M\subset S$ be a maximal spherical simplex containing $B$. 
Then there is a unique maximal spherical simplex $M'$ of $(W,S')$ 
and an element $u$ of $W$ such that 
$\langle M'\rangle = u\langle M\rangle u^{-1}$ by Proposition 4.15.  
By the Matching Theorem for systems of a finite Coxeter group applied 
to $(\langle M'\rangle,uMu^{-1})$ and $(\langle M'\rangle, M')$,  
there is a base $B'$ of $(\langle M'\rangle, M')$ such that 
$$[\langle B'\rangle,\langle B'\rangle]=
[\langle uBu^{-1}\rangle,\langle uBu^{-1}\rangle]=u[\langle B\rangle,\langle B\rangle]u^{-1}.$$ 
Moreover, $B'$ satisfies conditions 2 and 3, and so $|B|=|B'|$.

Let $C'$ be a base of $(W,S')$ that contains $B'$. 
Then by the above argument, there is a $C\subset S$ and a $v\in W$ such that 
$\langle C\rangle$ is a finite irreducible subgroup of $(W,S)$, and $|C| = |C'|$, and 
$[\langle C\rangle,\langle C\rangle]=v[\langle C'\rangle,\langle C'\rangle]v^{-1}.$ 
Then we have
$vu[\langle B\rangle,\langle B\rangle]u^{-1}v^{-1}\subset [\langle C\rangle,\langle C\rangle].$
By Lemma 4.18, there is a $w\in W$ such that $wBw^{-1}\subset C$. 
As $B$ is a base of $(W,S)$, we have 
that $wBw^{-1} = C = B$ by Lemma 4.16. Therefore $B'=C'$ 
and $B'$ is a base of $(W,S')$. 

Suppose $D'\subset S'$ is irreducible and $x\in W$ such that 
$[\langle D'\rangle,\langle D'\rangle]=x[\langle B\rangle,\langle B\rangle]x^{-1}.$
Then
$xu^{-1}[\langle B'\rangle,\langle B'\rangle]ux^{-1}=[\langle D'\rangle,\langle D'\rangle].$ 
By Lemma 4.18, there is a $y\in W$ such that $yB'y^{-1}\subset D'$. 
As $B'$ is a base of $(W,S')$ and $D'$ is irreducible, 
$yB'y^{-1} = D'= B'$ by Lemma 4.16. 
Thus $B'$ is unique. 
\end{proof}

\section{Nonisomorphic Basic Matching} 

Let $W$ be a finitely generated Coxeter group with two sets of Coxeter generators $S$ and $S'$. 
A base $B$ of $(W,S)$ is said to {\it match} a base $B'$ of $(W,S')$ 
if $[\langle B\rangle,\langle B\rangle]$ is conjugate to 
$[\langle B'\rangle,\langle B'\rangle]$ in $W$.  
In this section, we determine some necessary and some sufficient conditions 
for a base $B \subset S$ to match a base $B'\subset S'$ of a different type.

\begin{proposition} {\rm (Solomon \cite{Solomon}, Lemma 2)} 
If $A, B\subset S$ and 
$u$ is the shortest element of $\langle A\rangle u \langle B\rangle$, then 
$$\langle A\rangle\cap u\langle B\rangle u^{-1} = \langle C\rangle
\quad \hbox{where}\ \ C=A\cap uBu^{-1}.$$ 
\end{proposition} 

The next proposition follows from Proposition 5.1. 

\begin{proposition} 
Let $A, B\subset S$ and $w\in W$. 
Write $w = xuy$ with $x\in \langle A\rangle$, 
$y \in \langle B\rangle$, and $u$ the shortest element of $\langle A\rangle w\langle B\rangle$. 
Then 
$$\langle A\rangle\cap w\langle B\rangle w^{-1} = x\langle C\rangle x^{-1}
\quad \hbox{where}\ \ C=A\cap uBu^{-1}.$$ 
\end{proposition}

\begin{lemma} 
Suppose $B$ is a base of $(W,S)$ of type ${\bf C}_{2q+1}$ 
that matches a base $B'$ of $(W,S')$ of type ${\bf B}_{2q+1}$ for some $q\geq 1$. 
Let $a,b,c$ be the elements of $B$ such that $m(a,b)=4$ and $m(b,c)=3$, 
and let $A\subset S$ such that $a\in A$. 
If $\langle A\rangle$ is conjugate to $\langle A'\rangle$ 
for some $A'\subset S'$, then $B\subset A$. 
\end{lemma}
\begin{proof}
Let $M\subset S$ be a maximal spherical simplex containing $B$. 
Then there is a unique maximal spherical simplex $M'\subset S'$ 
such that $\langle M\rangle$ is conjugate to $\langle M'\rangle$ 
by Proposition 4.15. 
By conjugating $S'$, we may assume that $\langle M\rangle = \langle M'\rangle$.  
Then $M'$ contains $B'$ by the Basic Matching Theorem. 
Let $w$ be an element of $W$ such that $\langle A\rangle = w\langle A'\rangle w^{-1}$. 
By Proposition 5.2, there is an element $x$ of $\langle M'\rangle$ 
and a subset $C$ of $M'$ such that 
$$\langle M\cap A\rangle =\langle M\rangle\cap\langle A\rangle 
= \langle M'\rangle \cap w\langle A'\rangle w^{-1} = x\langle C\rangle x^{-1}. $$
Hence, we may assume that $W$ is finite by restricting to $\langle M\rangle$. 
Furthermore, by conjugating $S'$, we may assume that $\langle A\rangle =\langle A'\rangle$. 

Let $C$ be a base of $(W,S)$ other than $B$. 
Then each element of $C$ commutes with each element of $B$, and so   
$\langle B\rangle$ injects into the quotient of $W$ 
by the commutator subgroup of $\langle C\rangle$. 
Hence, by Theorem 3.9, 
we may assume that $W$ is the direct product of $\langle B\rangle$ and copies of ${\bf A}_1$. 
The center $Z$ of $W$ is generated by $S-B$ and 
the longest element $z$ of $\langle B\rangle$. 
The center $Z$ is also generated by $S'-B'$. 
Let $K$ be the kernel of the homomorphism of $Z$ to $\{\pm 1\}$ 
that maps $S-B$ to $1$ and $z$ to $-1$. Then $W/K$ is a Coxeter group. 
P-diagrams for $W/K$ are obtained from the P-diagram of $(W,S)$ 
by removing the vertices $S-B$ and from the P-diagram of $(W,S')$ 
by removing the vertices in $(S'-B')\cap K$ and 
identifying the remaining vertices of $S'-B'$ to a single vertex.  
By passing to the quotient $W/K$, 
we may assume that $S = B$ and $S'=B'\cup\{z'\}$ 
and $z'$ commutes with each element of $B'$. 
Then $z=z'$, since $\langle z\rangle =Z(W)=\langle z'\rangle$. 

Now as $a\in A$ and $\langle A\rangle = \langle A'\rangle$, 
we have that $a\in\langle A'\rangle$. 
The element $a$ is represented by the matrix ${\rm diag}(1,\ldots,1,-1)$ in ${\bf C}_{2q+1}$. 
Observe that
$${\rm diag}(1,\ldots,1,-1) = {\rm diag}(-1,\ldots,-1,-1){\rm diag}(-1,\ldots,-1,1).$$ 
The matrix $-I$ represents $z$. 
The matrix ${\rm diag}(-1,\ldots,-1,1)$ is the longest element of 
${\bf B}_{2q+1}$ which is in $[{\bf C}_{2q+1},{\bf C}_{2q+1}]$. 
Hence ${\rm diag}(-1,\ldots,-1,1)$ represents an element $\ell$ 
of $[\langle B\rangle, \langle B\rangle]$, with $a=\ell z$. 
As $[\langle B\rangle, \langle B\rangle]=[\langle B'\rangle, \langle B'\rangle]$, 
we have $\ell\in\langle B'\rangle$. 
Hence every reduced word in the generators $S'=B'\cup\{z\}$
representing $a = \ell z$ involves $z$ by Prop. 7 on p. 19 of Bourbaki\cite{Bourbaki}. 
Therefore $z\in A'$, 
since $A'\subset B'\cup\{z\}$ and $a\in\langle A'\rangle$. 
Hence $z\in \langle A\rangle$.  
As $z$ involves all the elements of $B$, we deduce that $B\subset A$. 
\end{proof}

\begin{theorem} 
Suppose $B$ is a base of $(W,S)$ of type ${\bf C}_{2q+1}$ 
that matches a base $B'$ of $(W,S')$ of type ${\bf B}_{2q+1}$ for some $q\geq 1$. 
Let $a,b,c$ be the elements of $B$ such that $m(a,b)=4$ and $m(b,c)=3$. 
If $s\in S-B$ and $m(s,a)<\infty$, then $m(s,t) = 2$ for all $t$ in $B$. 
\end{theorem}
\begin{proof}
Let $A\subset S$ be a maximal spherical simplex containing $\{a,s\}$.  
Then there is a maximal spherical simplex $A'\subset S'$ such that 
$\langle A\rangle$ is conjugate to $\langle A'\rangle$. 
Hence $B\subset A$ by Lemma 5.3. 
As $B$ is a base of $(\langle A\rangle,A)$, 
we deduce that $s$ commutes with each element of $B$. 
\end{proof}
\begin{lemma} 
Let $w=s_1\cdots s_n$ be a reduced word in $(W,S)$ 
and let $s\in S$ such that $s\neq s_i$ for each $i=1,\ldots,n$. 
If $sw$ has finite order in $W$, 
then $m(s,s_i)<\infty$ for each $i=1,\ldots,n$.
\end{lemma}
\begin{proof}
On the contrary, suppose $m(s,s_i) = \infty$ for some $i$. 
We may assume $S=\{s,s_1,\ldots,s_n\}$. 
Then 
$$W=\langle s, s_1,\ldots,\hat s_i, \ldots,s_n\rangle 
\ast_{\textstyle \langle s_1,\ldots,\hat s_i,\ldots,s_n\rangle}\langle s_1,\ldots,s_n\rangle$$
is a free product with amalgamation decomposition. 
Observe that $(sw)^k= swsw\cdots sw$ is a normal form for $(sw)^k$ 
for each $k\geq 1$ with respect to the amalgamated product, 
and so $(sw)^k \neq 1$ for each $k\geq 1$.
\end{proof}

\begin{theorem} 
Let $B$ be a base of $(W,S)$ of type ${\bf C}_{2q+1}$ for some $q\geq 1$, and   
let $a,b,c$ be the elements of $B$ such that $m(a,b)=4$ and $m(b,c)=3$. 
Suppose that $m(s,t)=2$ for all $(s,t)\in (S-B)\times B$ such that $m(s,a)<\infty$. 
Let $d =aba$, and let $z$ be the longest element of $\langle B\rangle$. 
Let $S'=(S-\{a\})\cup\{d,z\}$ and $B'=(B-\{a\})\cup\{d\}$. 
Then $S'$ is a set of Coxeter generators for $W$ 
such that 
\begin{enumerate}
\item The set $B'$ is a base of $(W,S')$ of type ${\bf B}_{2q+1}$ that matches $B$,  
\item $m(z,t)=2$ for all $t\in B'$,
\item If $(s,t)\in (S-B)\times\{d,z\}$, then 
$m(s,t)<\infty$ if and only if $m(s,a)<\infty$, moreover  
if $m(s,t)<\infty$, then $m(s,t)=2$.
\end{enumerate}
\end{theorem}
\begin{proof}
Consider the Coxeter presentation 
$$W =\langle S \ |\ (st)^{m(s,t)}:\, s,t \in S\ \hbox{and}\ m(s,t)<\infty\rangle$$
Now $(\langle B'\rangle, B')$ is a finite Coxeter system of type ${\bf B}_{2q+1}$. 
Let $\ell$ be the longest element of $(\langle B'\rangle, B')$. 
Regard $\ell$ as a reduced word in the elements of $B'$. 
Add generators $d$ and $z$ and relations $d=aba$ and $z=a\ell$ 
to the above presentation for $W$. 
Now add the relators $(st)^{m(s,t)}$ for $(s,t)$ in  $\{d,z\}\times S'$ or in  $S'\times\{d,z\}$ 
where $m(s,t)$ is the order of $st$ in $W$ and $m(s,t)<\infty$. 
This includes all the relators of $(\langle B'\rangle, B')$. 
As $\langle z\rangle$ is the center of $\langle B\rangle$, we have 
that $m(z,t)=2$ for all $t$ in $B'$. 

Next delete the generator $a$ and the relation $z=a\ell$ and replace $a$ by 
$z\ell$ in the remaining relations. 
As $z$ commutes with each element of $B'$, we can replace the relation $d=z\ell bz\ell$ by 
the relation $d=\ell b\ell$. 

The relators $(z\ell b)^4$ and $(bz\ell)^4$ can be replaced by 
$(\ell b)^4$ and $(b\ell)^4$ which in turn can be replaced by 
$(db)^2$ and $(bd)^2$ using the relation $d=\ell b\ell$. 
The relators $(db)^2$ and $(bd)^2$ are redundant and so we delete them. 
The relation $d=\ell b\ell$ is derivable from the relators 
of $(\langle B'\rangle, B')$ and so we delete it. 
The relators $(z\ell s)^2$ and $(sz\ell)^2$ for $s\in B-\{a,b\}$ 
can be replaced by $(\ell s)^2$ and $(s\ell)^2$. 
The relators $(\ell s)^2$ and $(s\ell)^2$ are derivable from the relators 
of $(\langle B'\rangle, B')$ and so we delete them.

Suppose $s\in S-B$ and $m(s,a) <\infty$. 
Then $m(s,t) = 2$ for all $t\in B$ by hypothesis. 
Hence $m(s,t) = 2$ for all $t\in B'\cup\{z\}$. 
Now the relators $(z\ell s)^2$ and $(sz\ell)^2$ can be replaced by 
$(\ell s)^2$ and $(s\ell)^2$. The relators $(\ell s)^2$ and $(s\ell)^2$ 
are derivable from the relators $(st)^2$ for $t\in B'$ 
and the relation $\ell^2=1$.
Hence we may delete the relators $(\ell s)^2$ and $(s\ell)^2$. 
This leaves the Coxeter presentation 
$$W =\langle S' \ |\ (st)^{m(s,t)}:\, s,t \in S'\ \hbox{and}\ m(s,t)<\infty\rangle$$
Thus $S'$ is a set of Coxeter generators for $W$. 

Statement 3 follows from Theorem 5.4, Lemma 5.5,  
and the hypothesis that $m(s,t)=2$ for all $(s,t)\in (S-B)\times B$ such that $m(s,a)<\infty$. 
The set $B'$ is a base of $(W,S')$, since $B$ is a base of $(W,S)$ 
and if $s\in S'-B'$ and $m(s,d) < \infty$, then $m(s,d)=2$. 
The base $B$ matches the base $B'$, 
since $[\langle B\rangle,\langle B\rangle] = [\langle B'\rangle,\langle B'\rangle]$. 
\end{proof}
 
The next theorem follows from Theorems 5.4 and 5.6. 

\begin{theorem} 
Let $(W,S)$ be a Coxeter system of finite rank.  
Let $B$ be a base of $(W,S)$ of type ${\bf C}_{2q+1}$ for some $q\geq 1$, and   
let $a,b,c$ be the elements of $B$ such that $m(a,b)=4$ and $m(b,c)=3$. 
Then $W$ has a set of Coxeter generators $S'$ such that $B$ matches a base $B'$ 
of $(W,S')$ of type ${\bf B}_{2q+1}$ if and only if 
$m(s,t)=2$ for all $(s,t)\in (S-B)\times B$ such that $m(s,a)<\infty$. 
\end{theorem} 

We next consider the analogue of Theorem 5.7 in the dihedral case. 

\begin{theorem} 
Let $(W,S)$ be a Coxeter system of finite rank, and   
let $B=\{a,b\}$ be a base of $(W,S)$ of type ${\bf D}_2(4q+2)$ for some $q\geq 1$. 
Then $W$ has a set of Coxeter generators $S'$ such that $B$ 
matches a base $B'$ of $(W,S')$ of type ${\bf D}_2(2q+1)$ if and only if 
either $v=a$ or $v=b$ has the property that 
if $s\in S-B$ and $m(s,v)<\infty$, then $m(s,a) = m(s,b) = 2$.  
\end{theorem}
\begin{proof}
Suppose that $W$ has a set of Coxeter generators $S'$ such that $B$ 
matches a base $B'$ of $(W,S')$ of type ${\bf D}_2(2q+1)$. 
Let $v=a$ or $b$ with the choice specified below. 
Suppose $s\in S-B$ and $m(s,v)<\infty$. 
Let $A\subset S$ be a maximal spherical simplex containing $\{s,v\}$. 
Then there is a maximal spherical simplex $A'\subset S'$ such that 
$\langle A\rangle$ is conjugate to $\langle A'\rangle$. 
We claim that $B\subset A$.  As in the proof of Lemma 5.3, we may assume  
that $\langle A\rangle = \langle A'\rangle$ and reduce $W$ so that $S=B$ and 
$S'=B'\cup\{z\}$ where $\langle z\rangle$ is the center of $\langle B\rangle$. 

Now $a$ and $b$ are not both in $\langle B'\rangle$. 
Choose $v$ so that $v$ is not in $\langle B'\rangle$. 
Then every reduced $S'$-word representing $v$ involves $z$. 
Now as $v\in A$, we have that $v\in \langle A'\rangle$. 
Therefore $z\in A'$. Hence $z\in \langle A\rangle$.  
Therefore $B\subset A$ as claimed. 
Now return to the original state of $W$. 
As $B$ is a base of $\langle A\rangle$, we have $m(s,a)=m(s,b)=2$. 
The converse follows from the next theorem. 
\end{proof}

\begin{theorem} 
Let $B=\{a,b\}$ be a base of $(W,S)$ of type ${\bf D}_2(4q+2)$ for some $q\geq 1$. 
Suppose that if $s\in S-B$ and $m(s,a)<\infty$, then $m(s,a) = m(s,b) = 2$. 
Let $c=aba$ and let $z$ be the longest element of $\langle B\rangle$. 
Let $S'=(S-\{a\})\cup\{c,z\}$ and $B'=\{b,c\}$. 
Then $S'$ is a set of Coxeter generators of $W$ such that 
\begin{enumerate}
\item The set $B'$ is a base of $(W,S')$ of type ${\bf D}_2(2q+1)$ that matches $B$,  
\item $m(z,b)=m(z,c) = 2$, 
\item if $(s,t)\in (S-B)\times \{c,z\}$, then $m(s,t)<\infty$ if and only if $m(s,a)<\infty$,  
moreover if $m(s,t)<\infty$, then $m(s,t) = 2$. 
\end{enumerate}
\end{theorem}
\begin{proof}
Consider the Coxeter presentation 
$$W =\langle S \ |\ (st)^{m(s,t)}:\, s,t \in S\ \hbox{and}\ m(s,t)<\infty\rangle$$
Now $(\langle B'\rangle, B')$ is a finite Coxeter system of type ${\bf D}_2(2q+1)$. 
Let $\ell$ be the longest element of $(\langle B'\rangle, B')$. 
Regard $\ell$ as the reduced word $(bc)^qb$ in the elements of $B'$. 
Add generators $c$ and $z$ and relations $c=aba$ and $z=a\ell$ 
to the above presentation for $W$. 
Now add the relators $(st)^{m(s,t)}$ for $(s,t)$ in  $\{c,z\}\times S'$ or in  $S'\times\{c,z\}$ 
where $m(s,t)$ is the order of $st$ in $W$ and $m(s,t)<\infty$. 
This includes all the relators of $(\langle B'\rangle, B')$. 
As $\langle z\rangle$ is the center of $\langle B\rangle$, we have 
that $m(z,b) = m(z,c) = 2$. 

Next delete the generator $a$ and the relation $z=a\ell$ and replace $a$ by 
$z\ell$ in the remaining relations. 
As $z$ commutes with each element of $B'$, we can replace the relation $c=z\ell bz\ell$ by 
the relation $c=\ell b\ell$. 

The relators $(z\ell b)^{2(2q+1)}$ and $(bz\ell)^{2(2q+1)}$ can be replaced by 
$(\ell b)^{2(2q+1)}$ and $(b\ell)^{2(2q+1)}$ which in turn can be replaced by 
$(cb)^{2q+1}$ and $(bc)^{2q+1}$ using the relation $c=\ell b\ell$. 
The relators $(cb)^{2q+1}$ and $(bc)^{2q+1}$ are redundant and so we delete them. 
The relation $c=\ell b\ell$ is derivable from the relators 
of $(\langle B'\rangle, B')$ and so we delete it. 

Suppose $s\in S-B$ and $m(s,a) <\infty$. 
Then $m(s,a) = m(s,b)= 2$ by hypothesis. 
Hence $m(s,t) = 2$ for $t\in\{b,c,z\}$. 
Now the relators $(z\ell s)^2$ and $(sz\ell)^2$ can be replaced by 
$(\ell s)^2$ and $(s\ell)^2$. The relators $(\ell s)^2$ and $(s\ell)^2$ 
are derivable from the relators $(st)^2$ for $t\in B'$ 
and the relation $\ell^2=1$.
Hence we may delete the relators $(\ell s)^2$ and $(s\ell)^2$. 
This leaves the Coxeter presentation 
$$W =\langle S' \ |\ (st)^{m(s,t)}:\, s,t \in S'\ \hbox{and}\ m(s,t)<\infty\rangle$$
Thus $S'$ is a set of Coxeter generators for $W$. 

Statement 3 follows from Lemma 5.5 and the hypothesis 
that $m(s,t) = m(s,b) = 2$ for all $s\in S-B$ such that $m(s,a)<\infty$. 
The set $B'$ is a base of $(W,S')$, since $B$ is a base of $(W,S)$ 
and if $s\in S'-B'$ and $m(s,c) < \infty$, then $m(s,c)=2$. 
The base $B$ matches the base $B'$, 
since $[\langle B\rangle,\langle B\rangle] = [\langle B'\rangle,\langle B'\rangle]$. 
\end{proof}

A group $G$ has property FA if for every tree on which $G$ acts without inversions, 
the set of fixed points of $G$ in the tree is nonempty. 
Let $(W,S)$ be a Coxeter system, and let $A\subset S$. 
We say that $\langle A\rangle$ is a {\it complete} visual subgroup of $(W,S)$ 
if $(\langle A\rangle,A)$ is a complete Coxeter system.

\begin{proposition} {\rm (Mihalik and Tschantz \cite{M-T})} 
Let $(W,S)$ be a Coxeter system of finite rank. 
The maximal {\rm FA} subgroups of $W$ are the conjugates of the maximal 
complete visual subgroups of $(W,S)$. 
\end{proposition}

\begin{lemma} 
Let $A,B\subset S$. If $\langle A\rangle$ is a maximal complete 
visual subgroup of $(W,S)$ and $\langle A\rangle$ and $\langle B\rangle$ are 
conjugate, then $A = B$. 
\end{lemma}
\begin{proof}
If $s\in S-A$, then the irreducible component of $\langle A\cup\{s\}\rangle$ 
containing $s$ is incomplete and therefore infinite, 
since $\langle A\rangle$ is a maximal complete visual subgroup of $(W,S)$.  
Hence no $s\in S-A$ is $A$-admissible, and so $A=B$ by Lemma 4.11.
\end{proof}

The next proposition follows form Proposition 5.10 and Lemma 5.11. 

\begin{proposition}  
Let $W$ be a finitely generated Coxeter group with two sets 
of Coxeter generators $S$ and $S'$, 
and let $M$ be a maximal simplex of $(W,S)$. 
Then there is a unique maximal simplex $M'$ of $(W,S')$ 
such that $\langle M\rangle$ and $\langle M'\rangle$ 
are conjugate in $W$. 
\end{proposition} 

Let $W$ be a finitely generated Coxeter group with two Coxeter systems $S$ and $S'$,  
and let $A$ be a subset of $S$. 
Let $\ov A$ be the intersection of all subsets $B$ of $S$ such that $B$ contains $A$ 
and $\langle B\rangle$ is conjugate to $\langle B'\rangle$ for some $B'\subset S'$. 
Then $\ov A$ is the smallest subset $B$ of $S$ such that $B$ contains $A$ 
and $\langle B\rangle$ is conjugate to $\langle B'\rangle$ for some $B'\subset S'$ 
by Prop. 5.2. If $A$ is a spherical simplex, then $\ov A$ is a spherical simplex,  
since for any maximal spherical simplex $M$ of $(W,S)$ that contains $A$,  
there exists $M'\subset S'$ such that 
$\langle M\rangle$ is conjugate to $\langle M'\rangle$ by Prop. 4.15.

\begin{theorem} 
Suppose $B =\{x,y\}$ is a base of $(W,S)$ of type ${\bf D}_2(2q+1)$ 
that matches a base $B'$ of $(W,S')$ of type ${\bf D}_2(4q+2)$ for some $q\geq 1$. 
Then there exists $r\in \ov B - B$ such that $\ov{\{r\}}=\ov B$.  
Moreover, if $s\in S-B$ and $m(s,x), m(s,y) <\infty$, then $m(s,x)=m(s,y)=2$.  
\end{theorem}
\begin{proof}
Let $C=\ov B$. 
Then $C$ is a spherical simplex of $(W,S)$ 
and $\langle C\rangle$ is conjugate to $\langle C'\rangle$ for some $C'\subset S'$. 
By conjugating $S'$, we may assume that $\langle C\rangle = \langle C'\rangle$. 
Then $C'$ contains $B'$ by the Basic Matching Theorem. 
Hence $B$ is a proper subset of $C$, since otherwise 
$\langle B'\rangle \subset \langle C'\rangle =\langle C\rangle = \langle B\rangle$ 
which is not the case, since $|\langle B\rangle|< |\langle B'\rangle|$. 

Let $r$ be an element of $C-B$ which will be specified below. 
As $r\in \ov B$, we have that $\ov{\{r\}}\subset \ov B$. 
Let $A= \ov{\{r\}}$. 
We claim that $B\subset A$. 
Now $\langle A\rangle$ is conjugate in $\langle C'\rangle$
to $\langle A'\rangle$ for some $A' \subset C'$ by Prop. 5.2. 
Hence we may assume that $S=C$. Then $W$ is a finite group. 
By quotienting out the commutator subgroups of all the bases of $(W,C)$ 
other than $B$, we may assume that $W$ is the direct product of $\langle B\rangle$ 
and copies of ${\bf A}_1$. 

The center $Z$ of $W$ is generated by $C-B$. 
The center $Z$ is also generated by $C'-B'$ and $z'$. 
Let $K$ be the kernel of the homomorphism $\phi:Z\to \{\pm 1\}$ 
that maps $C'-B'$ to $1$ and $z'$ to $-1$. 
Choose $r\in C-B$ so that $\phi(r)=-1$. 

By quotienting out $K$, we may assume that $C=B\cup\{r\}$ and $C'=B'$. 
Then $r=z'$, and so $z'\in \langle A\rangle$.  
Hence $z'\in \langle A'\rangle$. 
Therefore $A'=B'$, and so $B\subset A$ as claimed. 
Now return to the original state of $W$. 
As $B\subset A$, we have $\ov B \subset \ov{\{r\}}$. 
Thus $\ov{\{r\}}=\ov B$.  

Suppose $s\in S-B$ with $m(s,x), m(s,y) <\infty$. 
Let $M\subset S$ be a maximal simplex containing $\{s,x,y\}$. 
Then there is a maximal simplex $M'\subset S'$ such 
that $\langle M\rangle$ is conjugate to $\langle M'\rangle$ 
by Prop. 5.12. 
By conjugating $S'$, we may assume that $\langle M\rangle=\langle M'\rangle$. 
Then $M'$ contains $B'$ and $[B,B]$ is conjugate to $[B',B']$ in $\langle M'\rangle$ 
by the Basic Matching Theorem. 

Let $B'=\{a,b\}$. Then $m(s',a)=m(s',b)=2$ for all $s'\in M'-B'$ by Theorem 5.8. 
Hence $B'$ is an irreducible component of $M'$. 
Therefore $[B',B']$ is a normal subgroup of $\langle M'\rangle$. 
Hence $[B,B]$ is a normal subgroup of $\langle M\rangle$. 
As $\langle xy\rangle = [B,B]$, 
we have that $\langle xy\rangle$ is a normal subgroup of $\langle M\rangle$. 
Therefore $s\{x,y\}s =\{x,y\}$ by Lemma 4.17,   
and $sxs=x$ and $sys=y$ by the deletion condition. 
\end{proof}

\begin{lemma} 
Let $\phi:{\bf B}_n \to {\bf C}_n$ be a monomorphism with $n$ odd and $n\geq 3$. 
Then $\phi$ maps $b_{n-1}b_n$ to a conjugate of $(c_{n-1}c_n)^2$ in ${\bf C}_n$. 
\end{lemma}
\begin{proof}
Now $\phi({\bf B}_n)$ does not contain the center of ${\bf C}_n$, 
since $Z({\bf B}_n) = \{1\}$. 
Therefore either $\phi({\bf B}_n)={\bf B}_n$ or $\phi({\bf B}_n)=\theta({\bf B}_n)$ 
where $\theta$ is the automorphism of ${\bf C}_n$ defined by 
$\theta(c_i) = -c_i$, for $i=1,\ldots,n-1$ and $\theta(c_n)=c_n$. 
Now $\theta$ restricts to the identity on $[{\bf C}_n,{\bf C}_n]$, 
and so by composing $\phi$ with $\theta$ in the latter case, 
we may assume that $\phi({\bf B}_n)={\bf B}_n$. 
Now every automorphism of ${\bf B}_n$ is inner according to Franzsen [7]. 
Hence $\phi$ restricts to conjugation on $[{\bf B}_n,{\bf B}_n]$ 
by an element of ${\bf B}_n$. 
As $b_{n-1}b_n$ is in $[{\bf B}_n,{\bf B}_n]$ and $b_{n-1}b_n=(c_{n-1}c_n)^2$, 
we conclude that $\phi(b_{n-1}b_n)$ is conjugate to $(c_{n-1}c_n)^2$ in ${\bf C}_n$.
\end{proof}

\begin{theorem} 
Suppose $B$ is a base of $(W,S)$ of type ${\bf B}_{2q+1}$ 
that matches a base $B'$ of $(W,S')$ of type ${\bf C}_{2q+1}$ for some $q\geq 1$. 
Let $x,y$ be the split ends of the C-diagram of $(\langle B\rangle, B)$. 
Then there exists $r\in \ov B - B$ such that $\ov{\{r\}}=\ov B$.  
Moreover if $s\in S-B$ and $m(s,x), m(s,y) <\infty$, 
then $m(s,t)=2$ for all $t\in B$.  
\end{theorem}
\begin{proof}
The proof that there exists $r\in \ov B - B$ such that $\ov{\{r\}}=\ov B$ 
is the same as for Theorem 5.13.  
Suppose $s\in S-B$ with $m(s,x), m(s,y)<\infty$. 
Let $M\subset S$ be a maximal simplex containing $\{s,x,y\}$. 
Then there is a unique maximal simplex $M'\subset S'$ such 
that $\langle M\rangle$ is conjugate to $\langle M'\rangle$ 
by Prop. 5.12. 
By conjugating $S'$, we may assume that $\langle M\rangle=\langle M'\rangle$. 

Let $a, b, c$ be the elements of $B'$ such that $m(a,b)=4$ and $m(b,c)=3$. 
Now $xy$ is in $[\langle B\rangle,\langle B\rangle]$, 
and so $xy$ is conjugate to $(ab)^2$ by the 
Basic Matching Theorem and Lemma 5.14. 
Hence there is a $w\in W$ such that $w(ab)^2w^{-1}\in \langle M'\rangle$. 
Now $\langle (ab)^2\rangle =[\langle a,b\rangle,\langle a,b\rangle]$. 
Let $u$ be the shortest element of $\langle M'\rangle w\langle a,b\rangle$. 
Then $u\{a,b\}u^{-1}\subset M'$ by Lemma 4.18. 
As $m(a,b)=4$, we deduce that $\{a,b\}\subset M'$ by Lemma 4.11. 
Hence $B'\subset M'$ by Lemma 5.3.  Moreover $m(s',t') = 2$ 
for all $(s',t')\in (M'-B')\times B'$ by Theorem 5.4. 
Hence $B'$ is an irreducible component of $M'$. 
Therefore $[\langle B'\rangle,\langle B'\rangle]$ is a normal subgroup of $\langle M'\rangle$. 

Now $M$ contains $B$ and $[B,B]$ is conjugate to $[B',B']$ in $\langle M\rangle$ 
by the Basic Matching Theorem. 
Therefore $[\langle B\rangle,\langle B\rangle]=[\langle B'\rangle,\langle B'\rangle]$, 
since $[\langle B'\rangle,\langle B'\rangle]$ is a normal subgroup of $\langle M'\rangle$. 
Hence $[\langle B\rangle,\langle B\rangle]$ is a normal subgroup of $\langle M\rangle$. 
Then $sBs = B$ by Lemma 4.18, and $sts=t$ for all $t\in B$ by the deletion condition. 
\end{proof}

\section{Matching of Finite Irreducible Subgroups}  

As a reference for the automorphism groups of finite irreducible  
Coxeter groups, see Chapter 2 of Franzsen \cite{Franzsen} 
or \S6 of Franzsen and Howlett \cite{F-H}. 
An automorphism of a Coxeter system $(W,S)$ 
is called a {\it graph automorphism}. The graph automorphisms of $(W,S)$ 
correspond to the automorphisms of the P-diagram of $(W,S)$. 

\begin{lemma}  
Let $\alpha:{\bf B}_n \to {\bf B}_n$ be an automorphism. 
Then there is an inner automorphism $\iota$ of ${\bf B}_n$ and 
a graph automorphism $\gamma$ of ${\bf B}_n$ such that  
$\alpha|_{[{\bf B}_n,{\bf B}_n]}=\iota\gamma|_{[{\bf B}_n,{\bf B}_n]}$ with 
$\gamma$ the identity map if $n$ is odd. 
\end{lemma} 
\begin{proof}
If $n$ is odd, then every automorphism of ${\bf B}_n$ is inner. 
Assume that $n$ is even. 
Let $\psi$ be the automorphism of ${\bf B}_n$ defined by $\psi(w) = (-1)^{l(w)}w$. 
All the elements of $[{\bf B}_n,{\bf B}_n]$ have even length. 
Therefore $\psi$ restricts to the identity on $[{\bf B}_n,{\bf B}_n]$. 
Now there is an inner automorphism $\iota$ of ${\bf B}_n$ 
and a graph automorphism $\gamma$ of ${\bf B}_n$ such that 
$\alpha=\iota\gamma$ or $\alpha=\iota\gamma\psi$. 
Hence $\alpha|_{[{\bf B}_n,{\bf B}_n]}=\iota\gamma|_{[{\bf B}_n,{\bf B}_n]}$. 
\end{proof}

\begin{lemma}  
Let $\alpha:{\bf C}_n \to {\bf C}_n$ be an automorphism. 
Then there is an inner automorphism $\iota$ of ${\bf C}_n$ such that 
$\alpha|_{[{\bf C}_n,{\bf C}_n]}=\iota|_{[{\bf C}_n,{\bf C}_n]}$. 
\end{lemma} 
\begin{proof}
This is clear if $\alpha$ is inner, so suppose $\alpha$ is outer. 
Let $\theta$ be the automorphism of ${\bf C}_n$ defined 
by $\theta(c_i)=-c_i$, for $i=1,\ldots,n-1$, and $\theta(c_n)=c_n$. 
Then $\theta$ restricts to the identity on $[{\bf C}_n,{\bf C}_n]$. 
If $n$ is odd, then there is an inner automorphism $\iota$ of ${\bf C}_n$ 
such that $\alpha = \iota\theta$. 
Hence $\alpha|_{[{\bf C}_n,{\bf C}_n]}=\iota|_{[{\bf C}_n,{\bf C}_n]}$. 

Suppose now that $n$ is even. 
Let $\psi$ be the automorphism of ${\bf C}_n$ defined by $\psi(w) = (-1)^{l(w)}w$. 
All the elements of $[{\bf C}_n,{\bf C}_n]$ have even length. 
Therefore $\psi$ restricts to the identity on $[{\bf C}_n,{\bf C}_n]$. 
Now there is an inner automorphism $\iota$ of ${\bf C}_n$ such that either 
$\alpha = \iota\theta$, $\iota\psi$, or $\iota\theta\psi$. 
Hence $\alpha|_{[{\bf C}_n,{\bf C}_n]}=\iota|_{[{\bf C}_n,{\bf C}_n]}$. 
\end{proof}

\begin{lemma}  
Let $\alpha:{\bf F}_4 \to {\bf F}_4$ be an automorphism. 
Then there is an inner automorphism $\,\iota$ of ${\bf F}_4$ and 
a graph automorphism $\gamma$ of ${\bf F}_4$ such that  
$\alpha|_{[{\bf F}_4,{\bf F}_4]}=\iota\gamma|_{[{\bf F}_4,{\bf F}_4]}$. 
\end{lemma} 
\begin{proof}
Let $f_1,f_2,f_3,f_4$ be the Coxeter generators of ${\bf F}_4$ 
with $m(f_1,f_2) = 3$, $m(f_2,f_3)=4$, and $m(f_3,f_4) = 3$. 
Let $\psi_\ell$ be the automorphism of ${\bf F}_4$ defined by 
$\psi_\ell(f_i) = -f_i$ for $i=1,2$ and $\psi_\ell(f_i) = f_i$ for $i=3,4$. 
Then $\psi_\ell$ restricts to the identity on $[{\bf F}_4,{\bf F}_4]$. 
Now there is an inner automorphism $\iota$ of ${\bf F}_4$ 
and a graph automorphism $\gamma$ of ${\bf F}_4$ 
such that $\alpha = \iota\gamma$ or $\alpha=\iota\gamma\psi_\ell$. 
Hence $\alpha|_{[{\bf F}_4,{\bf F}_4]}=\iota\gamma|_{[{\bf F}_4,{\bf F}_4]}$. 
\end{proof}  

\begin{lemma}  
Let $\alpha:{\bf G}_4 \to {\bf G}_4$ be an automorphism. 
Then there is a reflection preserving automorphism $\beta$ of ${\bf G}_4$ such that  
$\alpha|_{[{\bf G}_4,{\bf G}_4]}=\beta|_{[{\bf G}_4,{\bf G}_4]}$. 
\end{lemma} 
\begin{proof}
Let $g_1,g_2,g_3,g_4$ be the Coxeter generators of ${\bf G}_4$ 
with $m(g_1,g_2) = 3$, $m(g_2,g_3)=3$, and $m(g_3,g_4) = 5$. 
According to Franzsen \cite{Franzsen}, the group ${\bf G}_4$ has an outer automorphism $\xi$ 
such that $\xi(g_i) = g_i$ for $i=1,2,3$ and $\xi(g_4)$ is conjugate to $g_4$. 
Let $\psi$ be the automorphism of ${\bf G}_4$ defined by $\psi(w) = (-1)^{l(w)}w$. 
All the elements of $[{\bf G}_4,{\bf G}_4]$ have even length. 
Therefore $\psi$ restricts to the identity on $[{\bf G}_4,{\bf G}_4]$. 
Now there is an inner automorphism $\iota$ of ${\bf G}_4$ 
such that $\alpha = \beta$ or $\beta\psi$ 
where $\beta=\iota$ or $\iota\xi$. 
Hence $\alpha|_{[{\bf G}_4,{\bf G}_4]}=\beta|_{[{\bf G}_4,{\bf G}_4]}$. 
\end{proof}

\begin{proposition} {\rm (Franzsen and Howlett \cite{F-H}, Prop. 32)} 
Let $(W,S)$ be a finite Coxeter system, and 
let $\alpha$ be an automorphism of $W$ that preserves reflections. 
Then $\alpha$ maps each visual subgroup 
of $(W,S)$ to a conjugate of a visual subgroup. 
\end{proposition}

Every automorphism of ${\bf A}_n$ is inner, except when $n = 5$. 
The group ${\rm Out}({\bf A}_5)$ has order two. 
The outer automorphisms of ${\bf A}_5$ behave badly with respect to visual subgroups 
because of the next proposition. 

\begin{proposition} {\rm (Franzsen and Howlett \cite{F-H}, Prop. 35)} 
Let $(W,S)$ be a finite Coxeter system of type ${\bf A}_5$, and 
let $\alpha$ be an automorphism of $W$.  
If there are proper subsets $A,B$ of $S$ and $w\in W$ such that 
$\alpha(\langle A\rangle) = w\langle B\rangle w^{-1}$, then $\alpha$ is inner. 
\end{proposition}

\begin{lemma} 
Let $(W,S)$ be a finite irreducible Coxeter system which is not of type ${\bf A}_5$. 
Let $\alpha$ be an automorphism of $W$, and let $A\subset S$. 
Then there is a $B\subset S$ such that $(\langle A\rangle, A) \cong (\langle B\rangle, B)$ 
and the group $\alpha([\langle A\rangle,\langle A\rangle])$ is conjugate to 
$[\langle B\rangle,\langle B\rangle]$ in $W$. 
\end{lemma}
\begin{proof}
If $(W,S)$ is of type ${\bf A}_n$,  
then $\alpha$ maps $\langle A\rangle$ to a conjugate of itself, 
since every automorphism of ${\bf A}_n$ is inner for all $n\neq 5$. 

Suppose $(W,S)$ is of type ${\bf B}_n$. 
By Lemma 6.1 there is an inner automorphism $\iota$ of $W$ and 
a graph automorphism $\gamma$ of $(W,S)$ such that $\alpha|_{[W,W]}=\iota\gamma|_{[W,W]}$. 
Let $B=\gamma(A)$. Then 
$\alpha([\langle A\rangle,\langle A\rangle])$ is conjugate to 
$[\langle B\rangle,\langle B\rangle]$. 

If $(W,S)$ is of type ${\bf C}_n$,  
then $\alpha$ maps $[\langle A\rangle,\langle A\rangle]$ to a conjugate of itself by Lemma 6.2. 
If $(W,S)$ is of type ${\bf D}_2(k)$, 
then $\alpha$ maps $\langle A\rangle$ to a conjugate of itself, 
since $\alpha$ preserves reflections. 
If $(W,S)$ is of type ${\bf E}_6$ or ${\bf E}_7$, 
then $\alpha$ maps $\langle A\rangle$ to a conjugate of itself, 
since every automorphism of ${\bf E}_6$ or ${\bf E}_7$ is inner. 

Suppose $(W,S)$ is of type ${\bf E}_8$. 
Let $\psi$ be the automorphism of $W$ defined by $\psi(w) = (w_0)^{l(w)}w$ 
where $w_0$ is the longest element of $(W,S)$. 
All the elements of $[W,W]$ have even length. 
Therefore $\psi$ restricts to the identity on $[W,W]$. 
Now there is an inner automorphism $\iota$ of $W$ such that 
$\alpha=\iota$ or $\iota\psi$. 
Hence $\alpha|_{[W,W]}=\iota|_{[W,W]}$. 
Therefore $\alpha$ maps $[\langle A\rangle,\langle A\rangle]$ to a conjugate of itself. 

Suppose $(W,S)$ is of type ${\bf F}_4$. 
By Lemma 6.3 there is an inner automorphism $\iota$ of $W$ and 
a graph automorphism $\gamma$ of $(W,S)$ such that $\alpha|_{[W,W]}=\iota\gamma|_{[W,W]}$. 
Let $B=\gamma(A)$. Then 
$\alpha([\langle A\rangle,\langle A\rangle])$ is conjugate to 
$[\langle B\rangle,\langle B\rangle]$. 
 
Suppose $(W,S)$ is of type ${\bf G}_3$. 
Then every automorphism of $W$ preserves reflections. 
Hence $\alpha$ maps $\langle A\rangle$ to a conjugate of itself by Prop. 6.5. 

Suppose $(W,S)$ is of type ${\bf G}_4$. 
By Lemma 6.4 there is 
a reflection preserving automorphism $\beta$ of $W$ such that 
$\alpha|_{[W,W]}=\beta|_{[W,W]}$.  
Therefore $\alpha$ maps $[\langle A\rangle,\langle A\rangle]$ to a conjugate of itself 
by Prop. 6.5. 
\end{proof}

The next proposition follows easily from Lemma 6.7. 

\begin{proposition} 
Let $(W,S)$ and $(W',S')$ be finite irreducible Coxeter systems  
which are not of type ${\bf A}_5$. 
Let $\alpha: W\to W'$ be an isomorphism, and let $A\subset S$. 
Then there is an $A'\subset S'$ such that $(\langle A\rangle, A) \cong (\langle A'\rangle, A')$ 
and the group $\alpha([\langle A\rangle,\langle A\rangle])$ is conjugate to 
$[\langle A'\rangle,\langle A\rangle]$ in $W'$. 
\end{proposition}

\begin{lemma} 
Let $n$ be odd with $n\geq 3$, and let $k$ be such that $3\leq k\leq n$. 
Identify ${\bf B}_k$ with $\langle b_{n-k+1},\ldots, b_n\rangle$ in ${\bf B}_n$ 
and ${\bf C}_k$ with $\langle c_{n-k+1},\ldots, c_n\rangle$ in ${\bf C}_n$. 
Let $\phi:{\bf B}_n\to {\bf C}_n$ be a monomorphism. Then $\phi$ maps 
$[{\bf B}_k,{\bf B}_k]$ to a conjugate of $[{\bf C}_k,{\bf C}_k]$ for each $k=3,\ldots,n$. 
\end{lemma}
\begin{proof} This follows from the proof of Lemma 5.14, since 
$[{\bf B}_k,{\bf B}_k]=[{\bf C}_k,{\bf C}_k]$ for each $k=3,\ldots,n$. 
\end{proof}

A {\it subbase} of a Coxeter system $(W,S)$ is a subset $A$ of $S$ 
such that $\langle A\rangle$ is a noncyclic, nonmaximal, finite, irreducible 
subgroup of $(W,S)$. 

\begin{theorem}  
{\rm (Subbase Matching Theorem)}
Let $W$ be a finitely generated Coxeter group with two sets of Coxeter generators $S$ and $S'$. 
Let $A$ be a subbase of $(W,S)$. Let $B$ be a base of $(W,S)$ containing $A$, 
and let $B'$ be the base of $(W,S')$ that matches $B$. 
Suppose that $B$ is not of type ${\bf A}_5$ and if $|\langle B\rangle|>|\langle B'\rangle|$, 
suppose that $A$ is not of type ${\bf C}_2$. 
Then $B'$ contains a subbase $A'$ of $(W,S')$ such that $[\langle A\rangle,\langle A\rangle]$ 
is conjugate to $[\langle A'\rangle,\langle A'\rangle]$ in $W$. Moreover
\begin{enumerate}
\item If $|\langle A\rangle|=|\langle A'\rangle|$, then $A$ and $A'$ have the same type. 
\item If $|\langle A\rangle|<|\langle A'\rangle|$, then $|\langle B\rangle|<|\langle B'\rangle|$ 
and $A$ is of type ${\bf B}_k$ and $A'$ is of type ${\bf C}_k$ for some $k\geq 3$. 
\item If $|\langle A\rangle|>|\langle A'\rangle|$, then $|\langle B\rangle|>|\langle B'\rangle|$ 
and $A$ is of type ${\bf C}_k$ and $A'$ is of type ${\bf B}_k$ for some $k\geq 3$. 
\end{enumerate}
Furthermore, if $|\langle B\rangle|>|\langle B'\rangle|$ and $A$ is of type ${\bf C}_2$, 
then $B'$ is of type ${\bf B}_{2q+1}$ for some $q\geq 1$ and 
$[\langle A\rangle,\langle A\rangle]$ is conjugate to $\langle xy\rangle$ in $W$ 
where $\{x,y\}$ is the set of split ends of the C-diagram of $(\langle B'\rangle,B')$. 
\end{theorem}
\begin{proof}
Suppose $|\langle B\rangle|=|\langle B'\rangle|$. 
By the Basic Matching Theorem, 
$B$ and $B'$ have the same type and 
there is an isomorphism $\phi:\langle B\rangle\to \langle B'\rangle$ 
that restricts to conjugation on $[\langle B\rangle,\langle B\rangle]$ by an element $u$ of $W$. 
By Prop. 6.8, there is a $A'\subset B'$ such that 
$(\langle A\rangle, A) \cong (\langle A'\rangle, A')$ 
and $\phi([\langle A\rangle,\langle A\rangle])$ is conjugate to 
$[\langle A'\rangle,\langle A'\rangle]$ by an element $v$ of $\langle B'\rangle$. 
Then $A'$ is a subbase of $(W,S')$ of the same type as $A$ 
and $vu[\langle A\rangle,\langle A\rangle]u^{-1}v^{-1}=[\langle A'\rangle,\langle A'\rangle]$. 

Now suppose $|\langle B\rangle|<|\langle B'\rangle|$. 
By the Basic Matching Theorem, 
$B$ is of type ${\bf B}_{2q+1}$ and $B'$ is of type ${\bf C}_{2q+1}$ for some $q\geq 1$ and 
there is an monomorphism $\phi:\langle B\rangle\to \langle B'\rangle$ 
that restricts to conjugation on $[\langle B\rangle,\langle B\rangle]$ by an element $u$ of $W$. 
Let $\alpha:(\langle B\rangle,B)\to {\bf B}_{2q+1}$ and 
$\beta:(\langle B'\rangle,B')\to {\bf C}_{2q+1}$ be isomorphisms of Coxeter systems. 
If $A$ is of type ${\bf A}_k$, we may assume, if necessary, 
by conjugating $\langle A\rangle$ by the longest element of $(\langle B\rangle, B)$,  
that $\alpha(A)\subset \langle b_1,\ldots, b_{2q}\rangle$. 
Now $b_i=c_i$ for $i=1,\ldots, 2q$. 
If $A$ is of type ${\bf A}_k$, let $A'=\beta^{-1}\alpha(A)$. 
Then $(\langle A\rangle,A)\cong(\langle A'\rangle, A')$. 
If $A$ is of type ${\bf B}_k$, let $\langle A'\rangle=\beta^{-1}({\bf C}_k)$ 
where ${\bf C}_k$ is as in Lemma 6.9 and $A'\subset S'$.  
Then $A'$ is of type ${\bf C}_k$. 

By the proof of Lemma 5.14, we deduce that 
$\beta\phi\alpha^{-1}:{\bf B}_{2q+1}\to {\bf C}_{2q+1}$ maps 
$[\langle\alpha(A)\rangle,\langle\alpha(A)\rangle]$ to 
$g[\langle\beta(A')\rangle,\langle\beta(A')\rangle]g^{-1}$ for some $g$ in ${\bf C}_{2q+1}$. 
Let $v=\beta^{-1}(g^{-1})$. 
Then $\phi([\langle A\rangle,\langle A\rangle])=v^{-1}[\langle A'\rangle,\langle A'\rangle]v$, 
and so 
$$vu[\langle A\rangle,\langle A\rangle]u^{-1}v^{-1}=[\langle A'\rangle,\langle A'\rangle].$$
The proof of the case $|\langle B\rangle|>|\langle B'\rangle|$ is the same as for the case 
$|\langle B\rangle|<|\langle B'\rangle|$ with the roles of $B$ and $B'$ reversed. 

Suppose $|\langle B\rangle|>|\langle B'\rangle|$ and $A$ is of type ${\bf C}_2$. 
By the Basic Matching Theorem, $B$ is of type ${\bf C}_{2q+1}$ and $B'$ 
is of type ${\bf B}_{2q+1}$ for some $q\geq 1$ 
and there is a monomorphism $\phi:\langle B'\rangle \to \langle B\rangle$ 
that restricts to conjugation on $[\langle B'\rangle,\langle B'\rangle]$ by 
an element of $W$. 
Let $\{x,y\}$ be the set of split ends of the C-diagram of $(\langle B'\rangle, B')$.  
Then $[\langle A\rangle,\langle A\rangle]$ is conjugate 
to $\langle xy\rangle$ by Lemma 5.14. 
\end{proof}

\begin{lemma} 
Let $(W,S)$ be a Coxeter system with $A,B \subset S$ such that 
$\langle A\rangle$ is finite and irreducible. 
If $\langle A\rangle$ is conjugate to $\langle B\rangle$ in $W$ 
and $A$ is neither of type ${\bf A}_n$, for some $n$,  
nor of type ${\bf B}_5$, then $A=B$. 
\end{lemma}
\begin{proof}
Suppose $s\in S-A$, with $m(s,a) > 2$ for some $a\in A$, and suppose $s$ is $A$-admissible. 
Then $K=A\cup\{s\}$ is irreducible. 
By Lemma 4.11, it suffices to show that $w_Ksw_K = s$. 
This is clear if $\langle w_K\rangle$ is the center of $\langle K\rangle$. 
Suppose that $Z(\langle K\rangle) = 1$. 
Now $K$ is not of type ${\bf A}_{n+1}$ nor of type ${\bf E}_6$, 
since $A$ is not of type ${\bf A}_n$ nor of type ${\bf B}_5$. 
Hence $K$ must be of type ${\bf B}_{2q+1}$ for some $q\geq 2$. 
Then $A$ is of type ${\bf B}_{2q}$ and $w_Ksw_K=s$. 
\end{proof}

\begin{theorem} 
Let $W$ be a finitely generated Coxeter group with two sets of Coxeter 
generators $S$ and $S'$. 
Let $A$ be a subbase of $(W,S)$, and let $A'\subset S'$.  
If $A'$ is irreducible and $[\langle A\rangle,\langle A\rangle]$ 
is conjugate to $[\langle A'\rangle,\langle A'\rangle]$ in $W$, 
then $A'$ is unique up to conjugation in $W$; 
moreover, if $A'$ is neither of type ${\bf A}_n$, for some $n$, nor of type ${\bf B}_5$, 
then $A'$ is unique. 
If $A$ is of type ${\bf C}_2$ and $A'=\{x,y\}$  
and $[\langle A\rangle,\langle A\rangle]$ is conjugate 
to $\langle xy\rangle$ in $W$, then $A'$ is unique. 
\end{theorem}
\begin{proof}
Suppose $A', A_1'$ are irreducible subsets of $S'$ and $[\langle A\rangle,\langle A\rangle]$ 
is conjugate to $[\langle A'\rangle,\langle A'\rangle]$ 
and to $[\langle A_1'\rangle,\langle A_1'\rangle]$. 
Then $[\langle A'\rangle,\langle A'\rangle]$ is conjugate to 
$[\langle A_1'\rangle,\langle A_1'\rangle]$. 
Hence $A'$ is conjugate to $A_1'$ by Lemma 4.18. 
If $A'$ is neither of type ${\bf A}_n$, for some $n$, nor of type ${\bf B}_5$, 
then $A'$ is unique by Lemma 6.11. 

Suppose $A$ is of type ${\bf C}_2$ and $A'=\{x,y\}\subset S'$ and 
$[\langle A\rangle,\langle A\rangle]$ is conjugate to $\langle xy\rangle$ in $W$. 
Then $m(x,y)=2$, since $[\langle A\rangle,\langle A\rangle]$ has order 2. 
Let $B$ be a base of $(W,S)$ containing $A$ and let $B'$ be the base of $(W,S')$ 
that matches $B$. Suppose $|\langle B\rangle| = |\langle B'\rangle|$. 
Then $B'$ contains a subbase $A'_1$ of type ${\bf C}_2$ 
such that $[\langle A\rangle,\langle A\rangle]$ is conjugate to 
$[\langle A_1'\rangle,\langle A_1'\rangle]$ in $W$ by Theorem 6.10. 
Then $xy$ is conjugate to an element of $\langle A_1'\rangle$. 
Hence $A'$ is conjugate to $A_1'$ by Lemma 4.17, 
which is a contradiction, since $A'$ and $A_1'$ have different types. 
Therefore $|\langle B\rangle| > |\langle B'\rangle|$ by Theorem 6.10. 

Now by the Basic Matching Theorem, 
$B$ is of type ${\bf C}_{2q+1}$ and $B'$ is of type ${\bf B}_{2q+1}$ for some $q\geq 1$. 
Let $E'=\{u,v\}$ be the set of split ends of the C-diagram of $(\langle B'\rangle, B')$. 
Then $[\langle A\rangle, \langle A\rangle]$ is conjugate to $\langle uv\rangle$ by Lemma 5.14. 
Hence $uv$ is conjugate to $xy$. 
Therefore $E'$ is conjugate to $A'$ by Lemma 4.17. 

Suppose $s\in S'-E'$ is $E'$-admissible. 
Then either $m(s,u)=m(s,v)=2$ or $m(s,u)=m(s,v)=3$ 
by Theorem 5.15. 
Let $K$ be the component of $E'\cup\{s\}$ containing $s$. 
Then $w_Ksw_K = s$. 
Therefore $E'=A'$ by Lemma 4.11. 
Thus $A'$ is unique. 
\end{proof}

\vspace{.2in}

\begin{theorem} {\rm (Edge Matching Theorem)} 
Let $W$ be a finitely generated Coxeter group with two sets 
of Coxeter generators $S$ and $S'$. 
Let $E=\{a,b\}$ be an edge of the P-diagram of $(W,S)$ 
with $m(a,b)\geq 4$. 
Then there is a unique edge $E'=\{x,y\}$ of the P-diagram of $(W,S')$ 
such that $[\langle E\rangle,\langle E\rangle]$ is conjugate in $W$ to either
$[\langle E'\rangle, \langle E'\rangle]$ or $\langle xy\rangle$. 
\end{theorem}
\begin{proof}
Assume first that $E$ is a base of $(W,S)$. 
Let $E'$ be the base of $(W,S')$ that matches $E$. 
Then $E'=\{x,y\}$ is an edge of the P-diagram $\Gamma'$ of $(W,S')$ 
such that $m(x,y)\geq 3$ by Theorem 4.19. 
To see that $E'$ is unique, 
suppose $E'_1=\{x_1,y_1\}$ is an edge of $\Gamma'$ such that 
$[\langle E\rangle, \langle E\rangle]$ is conjugate to either 
$[\langle E'_1\rangle, \langle E'_1\rangle]$ or $\langle x_1y_1\rangle$. 
Then $[\langle E'\rangle,\langle E'\rangle]$ is conjugate to either 
$[\langle E'_1\rangle, \langle E'_1\rangle]$ or $\langle x_1y_1\rangle$. 
Hence $E'$ is conjugate to $E'_1$ by Lemma 4.17 or Lemma 4.18. 
Therefore $E'=E'_1$ by Lemma 4.16. 

Assume now that $E$ is a subbase of $(W,S)$. 
Let $B$ be a base of $(W,S)$ containing $E$ and let $B'$ be the base of $(W,S')$ 
that matches $B$. 
Assume first that if $|\langle B\rangle|>|\langle B'\rangle|$, then $m(a,b)>4$. 
By Theorem 6.10, $B'$ contains a subbase $E'$ of $(W,S')$ of the same type as $E$ 
such that $[\langle E\rangle, \langle E\rangle]$ is conjugate to 
$[\langle E'\rangle, \langle E'\rangle]$. 
Suppose $E'_1=\{x_1,y_1\}$ is an edge of $\Gamma'$ such that 
$[\langle E\rangle, \langle E\rangle]$ is conjugate to either 
$[\langle E_1'\rangle, \langle E_1'\rangle]$ or $\langle x_1y_1\rangle$. 
As in the previous case, $E'$ is conjugate to $E'_1$. 
Therefore $E'=E'_1$ by Lemma 6.11. 

Now assume that $|\langle B\rangle|>|\langle B'\rangle|$ and $m(a,b)=4$. 
By Theorem 6.10, the base $B'$ is of type ${\bf B}_{2q+1}$ for some $q\geq 1$ and 
$[\langle E\rangle, \langle E\rangle]$ is conjugate to $\langle xy\rangle$ 
where $E'=\{x,y\}$ is the set of split ends of the C-diagram of $(\langle B'\rangle, B')$. 
Suppose $E'_1=\{x_1,y_1\}$ is an edge of $\Gamma'$ such that 
$[\langle E\rangle, \langle E\rangle]$ is conjugate to either 
$[\langle E'_1\rangle, \langle E'_1\rangle]$ or $\langle x_1y_1\rangle$. 
Then $\langle xy\rangle$ is conjugate to either 
$[\langle E'_1\rangle, \langle E'_1\rangle]$ or $\langle x_1y_1\rangle$. 
Hence $E'$ is conjugate to $E'_1$ by Lemma 4.17. 
Therefore $m(x_1,y_1)=2$ and $[\langle E\rangle, \langle E\rangle]$ is conjugate to 
$\langle x_1y_1\rangle$. Hence $E'=E'_1$ by Theorem 6.12.
\end{proof}

\section{Visual Graph of Groups Decompositions}  

Let $(W,S)$ be a Coxeter system of finite rank.  
Suppose that $S_1,S_2\subset S$, with $S=S_1\cup S_2$, and $S_0=S_1\cap S_2$ 
are such that there is no defining relator of $W$ (no edge of
the P-diagram) between an element of $S_1-S_0$ and $S_2-S_0$. 
Then we can write $W$ as a visual amalgamated product
$W=\langle S_1\rangle*_{\langle S_0\rangle}\langle S_2\rangle$. 
We say that $S_0$ {\it separates} $S$ if
$S_1-S_0\neq\emptyset$ and $S_2-S_0\neq\emptyset$.  
The amalgamated product decomposition of $W$ will be
nontrivial if and only if $S_0$ separates $S$. 
If $S_0$ separates $S$, we call the triple $(S_1,S_0,S_2)$ 
a {\it separation} of $S$. 
Note that $S_0$ separates $S$ if and only if $S_0$ 
separates $\Gamma(W,S)$, that is, there are $s_1,s_2$ in $S-S_0$ 
such that every path in $\Gamma(W,S)$ from $s_1$ to $s_2$ must pass 
through $S_0$. 

Let $\ell \in \langle S_0\rangle$ such that 
$\ell S_0\ell^{-1}=S_0$. 
By Lemma 4.5, we have $S_0=S_{\bullet}\cup (S_0-S_{\bullet})$
where $S_{\bullet}$ generates a finite group, each
element of $S_{\bullet}$ commutes with each element of $S_0-S_{\bullet}$, 
and $\ell$ is the longest element of $\langle S_{\bullet}\rangle$.  
The triple $(S_1,\ell, S_2)$ determines an elementary
twist of $(W,S)$ (or of its P-diagram) giving a new Coxeter
generating set $S_*=S_1\cup \ell S_2\ell^{-1}$ of $W$.

In application, it is simpler to consider a more general kind of twisting.  
Suppose $S_0$ and $\bar S_0\subset S_2$ generate conjugate subgroups of
$\langle S_2\rangle$.  Suppose $d\in \langle S_2\rangle$ is such
that $d\bar S_0d^{-1}=S_0$.  Then $S_1\cap dS_2d^{-1}=S_0$, since
$$S_0\subset S_1\cap dS_2d^{-1}\subset S_1\cap\langle S_2\rangle = S_0.$$
A {\it generalized twist} (or simply {\it twist}) of $(W,S)$ 
in this situation gives a new Coxeter generating
set $S_*=S_1\cup dS_2d^{-1}$ of $W$ and a new visual 
amalgamated product decomposition 
$W=\langle S_1\rangle*_{\langle S_0\rangle}\langle d S_2d^{-1}\rangle$

Elementary and generalized twists can be easily understood in terms 
of their effects on P-diagrams. 
The P-diagram of $(W,S)$ is the union of the 
P-diagrams for $\langle S_1\rangle$ and $\langle S_2\rangle$
overlapping in the P-diagram for $\langle S_0\rangle$. 
The P-diagram for $(W,S_*)$ is obtained from the P-diagram of $(W,S)$ 
by twisting the P-diagram of $\langle S_2\rangle$, that is, 
removing the P-diagram for $\langle S_2\rangle$, replacing it by the isomorphic 
P-diagram of $\langle dS_2d^{-1}\rangle$, and attaching it to the P-diagram
for $\langle S_1\rangle$ along $S_0 = d\bar S_0d^{-1}$.   
If $S_0 = \emptyset$, we call the twist {\it degenerate}. 
A degenerate twist does not change the isomorphism type of the P-diagram. 
This includes the case where
$S_1=S_0=\emptyset$, $S_2=S$, giving $S_*=dSd^{-1}$ the
conjugation of $S$ by an arbitrary $d\in W$.  
Any nondegenerate generalized twist of a Coxeter system $(W,S)$ can be realized 
by a finite sequence of elementary twists.

Let $\Lambda$ be a visual graph of groups decomposition of $(W,S)$. 
Then the graph of $\Lambda$ is a tree, since the abelianization of $W$ is finite. 
The graph of groups decomposition $\Lambda$ can
be understood as a visual amalgamated product in many ways, e.g.,
by taking some of the vertex and edge groups to be generated by
$S_1$, others to be generated by $S_2$, with the overlap being a
single edge group of $\Lambda$.  
Hence we will also speak of twisting a visual
graph of groups decomposition with respect to some such
partitioning of the graph of groups and some conjugating element.

A graph of groups decomposition is said to be {\it reduced} 
if no edge group is equal to an incident vertex group. 
Suppose $\Lambda$ is a reduced visual graph of groups
decomposition of a Coxeter system $(W,S)$ of finite rank.  
Suppose (for simplicity in this application) that no edge group of $\Lambda$ 
is a proper subgroup of another edge group of $\Lambda$.  
Construct another visual graph of groups decomposition 
(though not reduced) $\tilde\Lambda$ as follows. 
The vertices of $\tilde\Lambda$ are of two distinct types,  
v-vertices and e-vertices. 
The v-vertices correspond to the vertices of $\Lambda$, and the e-vertices 
correspond to the distinct edge groups of $\Lambda$.  
An edge of $\tilde\Lambda$ will connect vertices $p$ and $q$ 
if $p$ is a v-vertex and $q$ is an e-vertex, 
and $p$ corresponds to an endpoint of an edge of $\Lambda$ 
with edge group corresponding to $q$. 
The vertex group of a v-vertex $p$ of $\tilde\Lambda$ 
will be the vertex group for $p$ in $\Lambda$. 
The vertex group of an e-vertex $q$ of $\tilde\Lambda$ 
is the edge group of $\Lambda$ corresponding to $q$. 
Each edge of $\tilde\Lambda$,  say from
$p$ to $q$, will have edge group equal to the edge group of $\Lambda$ 
corresponding to the e-vertex $q$ of that edge.  
The maps of edge groups into vertex groups in $\tilde\Lambda$
will be inclusion maps (as in visual decompositions generally,
determined by which of the generators lie in each vertex and edge groups).  
Then by a series of reductions and expansions (inverse
reductions) we can get from $\tilde\Lambda$ to $\Lambda$ and we
see that they are both visual graph of group decompositions of $W$
(or by comparing the relations defining the fundamental groups of
$\Lambda$ and $\tilde\Lambda$). On the other hand, different
reduced visual graphs of groups $\Lambda$ and $\Lambda_2$ will
correspond to the same $\tilde\Lambda=\tilde\Lambda_2$ provided
they have the same vertex groups and edge groups, since the
inclusion of edge groups into vertex groups determine the edges of
$\tilde\Lambda$ when no edge group is a proper subgroup of another
edge group. The point here is that $\tilde\Lambda$ provides a way
of keeping track of which edge groups of $\Lambda$ are equal and
in which vertex groups without specifying what the subtree of
edges of $\Lambda$ with the same given edge group must look like,
(in essence, without specifying the order of the vertex groups
containing this edge group).  We think of $\tilde\Lambda$ as a
{\it flattened form} of $\Lambda$ making uniform the relationship
between vertex groups and different edge groups of $\Lambda$.

A particularly simple case is when the edge groups of $\Lambda$
are all equal.  Then $\tilde\Lambda$ has one v-vertex for each
vertex $p$ of $\Lambda$ and one e-vertex $q$ for the common edge
group $E$ with edges from $q$ to $p$ for each v-vertex $p$ and with edge
groups equal to $E$ and inclusion maps into the vertex groups. 
The fundamental group of $\Lambda$ is an amalgamated product 
of all the vertex groups of $\Lambda$ identifying the copies of the edge group in
each vertex group.

\section{The Decomposition Matching Theorem}  

If $U$ is a subgroup of $W$, write $U^*=\lbrace wUw^{-1}:\,w\in
W\rbrace$ for the set of all subgroups conjugate to $U$ in $W$.
Write $U^*\preceq V^*$ if for some $w\in W$, $U\subseteq wVw^{-1}$
(independent of the representatives for the conjugacy classes).
Clearly $\preceq$ is transitive and reflexive. 
Consider the conjugacy classes of a visual
subgroup $U$ and any subgroup $V$ of a Coxeter system $(W,S)$. If
$U^*\preceq V^*$ and $V^*\preceq U^*$ then $U^*=V^*$, since if
$U\subseteq wVw^{-1}\subseteq wzUz^{-1}w^{-1}$ then, since $U$ is
a visual subgroup, $wzUz^{-1}w^{-1}=U$ and $U$ and $V$ are
conjugate by Lemma 4.3.  Hence for the conjugacy
classes of visual subgroups, $\preceq$ is a partial order.  We say
that $J\subseteq S$ is a c-minimal separating subset of generators
if $\langle J\rangle^*$ is a $\preceq$-minimal element of the set
of conjugacy classes of subgroups generated by separating subsets of $S$.  
Assuming there are separating subsets of $S$,  
there are finitely many since $S$ is finite, and so
there are c-minimal separating subsets of $S$.

\begin{theorem}\label{Lambdadecomptheorem}  
Suppose $(W,S)$ and $(W,S')$ are two Coxeter systems
for the same finitely generated Coxeter group. 
If $(W,S)$ is complete, then $(W,S')$ is complete;  
otherwise, for any given nontrivial splitting $A*_CB$ of $W$, 
there exist $S_0\subseteq S$,
$S'_0\subseteq S'$, a visual graph of groups decomposition
$\Lambda$ for $(W,S)$, and a visual graph of groups decomposition
$\Lambda'$ for $(W,S')$ such that:
\begin{enumerate}
\item $S_0$ is a c-minimal separating subset of $S$, $S'_0$ is a
c-minimal separating subset of $S'$, with $\langle
S_0\rangle^*=\langle S'_0\rangle^*\preceq C^*$;

\item the edge groups of $\Lambda$ are conjugate to $\langle
S_0\rangle$, the edge groups of $\Lambda'$ are conjugate to
$\langle S'_0\rangle$ (and hence are conjugate and conjugate to a
subgroup of $C$); and

\item there is a 1-1 correspondence between the vertices of
$\Lambda$ and the vertices of $\Lambda'$ such that each vertex
group of $\Lambda$ is conjugate to the corresponding vertex group
of $\Lambda'$.
\end{enumerate}
\end{theorem}

\begin{proof}
If $(W,S)$ is complete, then $(W,S')$ is complete by Prop. 5.10. 
Suppose $(W,S)$ is incomplete. 
Given a nontrivial splitting $W=A*_CB$, there
is some visual splitting $W=A_1*_{C_1}B_1$, with respect to $S$,
with $C_1$ a subgroup of a conjugate of $C$ 
by the visual decomposition theorem and Remark 1 in \cite{M-T}. 
Consider the finite collection of conjugacy classes $\langle J\rangle^*$,
partially ordered by $\preceq$, for subsets $J\subseteq S$ such
that $\langle J\rangle^*\preceq C^*$ and there is a visual
splitting $W=A_2*_{\langle J\rangle}B_2$. Then there exists such a
$J$ with $\langle J\rangle^*$ minimal in this partial order.

Now starting with a splitting $W=A_2*_{\langle J\rangle}B_2$ and
working with respect to $S'$, as above, there is a $J'\subseteq
S'$ with $\langle J'\rangle^*\preceq\langle J\rangle^*$ and a
visual splitting $W=A_3*_{\langle J'\rangle}B_3$, with $\langle
J'\rangle^*$ $\preceq$-minimal for such splitting $S'$-visual
subgroups.

Working back again, from $W=A_3*_{\langle J'\rangle}B_3$ and
splitting visually with respect to $S$, there is a $J''\subseteq
S$ with $\langle J''\rangle^*\preceq\langle J'\rangle^*$ and an
$S$-visual splitting over $\langle J''\rangle$ with $\langle
J''\rangle^*$ $\preceq$-minimal.  Now $\langle
J''\rangle^*\preceq\langle J'\rangle^*\preceq\langle
J\rangle^*\preceq C^*$ but $J$ was taken so $\langle J\rangle^*$
was $\preceq$-minimal below $C^*$ having an $S$-visual splitting
over $\langle J\rangle$, hence $\langle J''\rangle^*=\langle
J\rangle^*$ (but not necessarily $\langle J''\rangle=\langle
J\rangle$), and so in fact $\langle J\rangle^*=\langle
J'\rangle^*$ and (1) holds with $S_0=J$ and $S'_0=J'$.

Since $S$ is finite and each vertex group of a reduced visual
graph of groups decomposition of $W$ is generated by a different
subset of the generators, there is an obvious limit to the number
of vertices in a reduced visual graph of groups decomposition of
$W$, and in some sense, the more vertices, the finer the graph of
groups decomposition. Take a reduced $S$-visual graph of groups
decomposition $\Lambda$ of $W$ such that every edge group is
conjugate to $\langle J\rangle$ and, among such, having a maximum
number of vertices. By the visual decomposition theorem, take
$\Lambda'$ a reduced $S'$-visual graph of groups decomposition
refining $\Lambda$, i.e., such that each vertex (edge) group of
$\Lambda'$ is a subgroup of a conjugate of a vertex (edge) group of
$\Lambda$. Similarly, take $\Lambda''$ a reduced $S$-visual graph
of groups decomposition of $W$ refining $\Lambda'$. 
The edge groups of $\Lambda''$ are equal to conjugates of the edge 
groups of $\Lambda$ by the c-minimality of $S_0$, and so are 
conjugate to the edge groups of $\Lambda'$, and so (2) holds. 
We postpone the proof of (3) until after the proof of Lemma 8.3. 

The following lemma characterizes the visual decomposition $\Lambda$. 

\begin{lemma} 
Suppose $(W,S)$ is a Coxeter system of finite rank 
and $J$ is a c-minimal separating subset of $S$. 
Let ${\mathcal E}$ be the set of separating subsets of $S$ 
that are conjugate to $J$ in $W$.  
Let ${\mathcal V}$ be the set of all maximal subsets of $S$ 
that are not separated by a set in ${\mathcal E}$. 
Suppose $\Lambda$ is a reduced visual graph of groups
decomposition of $(W,S)$ having edge groups generated by
conjugates of $J$ (and hence elements of ${\mathcal E}$) and among
such has a maximum number of vertices. Then all of the subgroups
generated by sets in ${\mathcal V}$ are the vertex groups of
$\Lambda$, and all of the subgroups generated by sets in
${\mathcal E}$ are the edge groups of $\Lambda$. 
\end{lemma}

\begin{proof}
All the visual conjugates of $\langle J\rangle$ are visual direct products $F\times G$ 
with conjugate finite factors $F$ and the same factor $G$ in common 
with all the visual conjugates of $\langle J\rangle$. 
If we split $W$ by a separating visual conjugate of $J$, 
each of the other visual conjugates of $J$
lies entirely in one of the factors of the free product with amalgamation, 
since the corresponding finite group $F$ 
lies in one factor and the group $G$ lies in each factor. 

Assume $\Lambda$ is a reduced visual graph of groups with edge
groups conjugate to $\langle J\rangle$ and among such having a maximal number of
vertices. Note that each edge group, and hence each vertex group,
contains the common subgroup $G$.  
The graph of $\Lambda$ is a tree, since the abelianization of $W$ is finite. 

Suppose $L$ is the set of generators of a vertex group $V$ of $\Lambda$. 
We claim that $L$ is not separated by a set in ${\mathcal E}$. 
On the contrary, suppose $L$ is separated by a set $K$ in ${\mathcal E}$, 
say $x$ and $y$ are in different components of the P-diagram of $\langle L-K\rangle$. 
We claim that $K\subseteq L$. On the contrary,  
suppose $K\not\subseteq L$. Then $L\cap K$ does not
separate $S$, by c-minimality of $J$, and so there is a path in the
P-diagram of $(W,S)$ from $x$ to $y$ that avoids $L\cap K$.  
Take a path from $x$ to $y$ which is in a union of as few vertex groups
of $\Lambda$ as possible. Let $V'$ be a vertex group of
$\Lambda$ containing a generator in this path not in $V$.  
Then the path passes through some edge group $E$ of $V$ at some first
point before $V'$ and must pass back through $E$ at some last point, 
since the graph of $\Lambda$ is a tree.  
Neither of these points is a generator of $G$
since these all lie in $L\cap K$.  Hence these points are
generators in the finite factor $F$ of $E$.  But the P-diagram of
$F$ is complete, and so there is a short circuit of the path going
from the first to the last point in $F$ avoiding $V'$.  
We conclude instead that the path hitting the fewest vertex groups of
$\Lambda$ is a path in $L-K$, contradicting the assumption that
$K$ separates $L$.  
Hence $K \subseteq L$ and there is a separation $(L_1,K,L_2)$ of $L$.  
Each edge group of $\Lambda$ incident to 
the vertex group $V = \langle L\rangle$ is contained in either the
subgroup generated by $L_1$ or by $L_2$, and so we can split $V$ 
into two vertices generated by $L_1$ and $L_2$, respectively, 
and joined by an edge group generated by $K$, 
with each component of the rest of $\Lambda$
attached to one or the other of the new vertex groups 
by an edge group of $\Lambda$. 
Neither of the new vertex groups equals an incident edge group $E$, 
since the finite Coxeter groups $E/G$ and $\langle K\rangle/G$ have the same rank.  
This gives a reduced visual graph of groups decomposition over 
separating conjugates of $J$ with more vertex groups, 
contradicting the maximality of the number of vertices in $\Lambda$.  
Hence $L$ cannot be separated by a set in ${\mathcal E}$ as claimed. 
Clearly, every subset of $S$ that contains $L$ properly  
is separated by the set of generators of some edge group of $\Lambda$ 
that is incident to $V$. 
Therefore $L$ is a maximal subset of $S$ 
that is not separated by a set in ${\mathcal E}$, 
and so $L\in {\mathcal V}$.

Now suppose $L\in{\mathcal V}$. 
We claim that $\langle L\rangle$ is a vertex group of $\Lambda$. 
Every element of $L$ is a generator of some vertex group of $\Lambda$.
Suppose $L'\subseteq L$ is a maximal subset of $L$ contained in
some vertex group of $\Lambda$. If $L-L'\neq\emptyset$, 
say $x\in L-L'$, then $L'$ and $x$ are not both contained in a vertex group
of $\Lambda$.  Take vertex groups $V$ and $V'$ of $\Lambda$,
with $x\in V$ and $L'\subseteq V'$, which are closest together
in the graph of $\Lambda$.  Let $E$ be an edge group of the path
between $V$ and $V'$. Then $E$ is generated by a visual conjugate $K$ of $J$
which separates the generators in $V-E$ from those in $V'-E$,
and so $K\in{\mathcal E}$.  Now $x\notin E$ otherwise $x$ would also
be in a vertex group closer to $V'$ on the path between $V$ and $V'$. 
Likewise, $L'\not\subseteq E$ or else $L'$ would be contained
in a vertex group closer to $V$ on a path between $V$ and $V'$. 
But then the P-diagram of $\langle L-K\rangle$ would have at least two components, 
one containing $x$ and one containing some element of $L'-K$.  
This contradicts the assumption that $L\in{\mathcal V}$.  
Instead all of $L$ must be contained in a vertex group $V$ of $\Lambda$. 
As the set of generators in $V$ is in ${\mathcal V}$, 
we have that $\langle L\rangle = V$.

Finally, suppose $K\in{\mathcal E}$.  
Then there is a separation $(S_1,K,S_2)$ of $S$. 
Each $L\in{\mathcal V}$ generates a vertex group of
$\Lambda$ but is not separated by $K$ by our previous argument, 
and so each $L\in{\mathcal V}$
is contained in either $S_1$ or $S_2$.  Pick vertex groups $V_1$
and $V_2$ as close together in $\Lambda$ as possible such that
$V_1$ is generated by a subset of $S_1$ and $V_2$ is generated by
a subset of $S_2$.  Then $V_1$, and $V_2$ are adjacent since every
vertex group in a path between these is generated by a subset of
either $S_1$ or $S_2$.  Now $V_1\cap V_2$ is an edge group $E$ of
$\Lambda$ which is generated by a subset of $K$ but not by a proper
subset of $K$ by the c-minimality of $J$, and so $E=\langle K\rangle$. 
\end{proof}

The next lemma explains the relationship between the 
visual decompositions $\Lambda$ and $\Lambda''$ of $(W,S)$. 

\begin{lemma}\label{LambdaLambdapplemma} 
Suppose $\Lambda$ and $\Lambda''$ are reduced visual graph of
groups decompositions of a Coxeter system $(W,S)$ of finite rank.  
Suppose the edge groups of $\Lambda$ are generated by
conjugates of a c-minimal separating subset $J$ of $S$, and, among
visual decompositions with this same conjugacy class of edge
groups, $\Lambda$ has a maximum number of vertex groups.  
Suppose each vertex and edge group of $\Lambda''$ is a subgroup of a
conjugate of a vertex or edge group of $\Lambda$, respectively.
Then the vertex and edge groups of $\Lambda$ are equal to the
vertex and edge groups of $\Lambda''$, respectively, that is,
$\tilde\Lambda=\tilde\Lambda''$.
\end{lemma}

\begin{proof}
By the last lemma, the vertex groups of $\Lambda$ are determined
from the set of all separating sets of generators 
that are conjugate to $J$. Each edge group of $\Lambda''$ is generated
by a separating subset of $S$ and is contained in a conjugate of a
$\langle J\rangle$, and so, by the c-minimality of $J$, must be a
conjugate of $\langle J\rangle$ and an edge group also of $\Lambda$. 

Let $T$ be the Bass-Serre tree with standard transversal $T_\ast$, 
corresponding to the graph of groups $\Lambda$,  
i.e., the vertices of $T$ are the cosets of each vertex group of $\Lambda$ and $T_\ast$ 
consists of the cosets of each vertex group that contain the identity.  
A vertex group $G$ of $\Lambda''$ stabilizes a vertex $V$ of $T$, 
since $G$ is a subgroup of a conjugate of a vertex group of $\Lambda$.  
But each generator of $G$ also stabilizes a vertex of $T_\ast$ 
and the geodesic path from that vertex of $T_\ast$ to $V$.  
Hence $G$ also stabilizes the vertex of $T_\ast$ nearest to $V$.  
Thus each vertex group of $\Lambda''$ is actually a subgroup 
of a vertex group of $\Lambda$.

As the vertex groups of $\Lambda$ are proper subgroups,
$\Lambda''$ has at least two vertices, and each vertex group of
$\Lambda''$ contains an edge group of $\Lambda''$, 
which is a conjugate of $\langle J\rangle$, as a proper subgroup.  
Hence no vertex group of $\Lambda''$ is contained in an edge group of $\Lambda$, 
since all the visual conjugates of $\langle J\rangle$ have the same rank. 
Consequently, each vertex group of $\Lambda''$ can be contained in
only one vertex group of $\Lambda$, otherwise 
a vertex group of $\Lambda''$ would be contained in the intersection of vertex
groups for two different vertices of $\Lambda$ and so would be contained in each edge
group for edges of $\Lambda$ in the geodesic path between these vertices, 
which is not the case. 

Summarizing, for each vertex $U$ of $\Lambda''$, there exists a
unique vertex $f(U)$ of $\Lambda$ such that the vertex group
$\Lambda''(U)$ of $\Lambda''$ at $U$ is a subgroup of the vertex
group $\Lambda(f(U))$ of $\Lambda$ at $f(U)$.  We claim that for
each vertex $V$ of $\Lambda$, the vertex group $\Lambda(V)$ is
generated by the vertex groups of $\Lambda''$ for vertices in
$f^{-1}(V)$.  In particular, there will be at least one vertex of
$\Lambda''$ in $f^{-1}(V)$, and so at least as many vertices in
$\Lambda''$ as in $\Lambda$.  But $\Lambda$ has a maximal number
of vertices for visual reduced graph of groups decompositions of
$(W,S)$ with edge groups that are conjugates of $\langle J\rangle$, so
$\Lambda''$, which also satisfies these conditions, has no more
vertices than $\Lambda$.  Hence $\Lambda$ and $\Lambda''$ have the
same number of vertices; moreover, for each vertex $V$ of $\Lambda$,
we conclude that $f^{-1}(V)$ is a unique vertex of $\Lambda''$, and the vertex
groups of these vertices in $\Lambda$ and $\Lambda''$ must be equal.  
Hence the vertex and edge groups of $\Lambda$ are the same
as the vertex and edge groups of $\Lambda''$, respectively, and so 
$\tilde \Lambda=\tilde\Lambda''$. 

To establish the claim that each vertex group $\Lambda(V)$ of
$\Lambda$ is generated by the vertex groups of $\Lambda''$ that it
contains, we will show that each edge group of $\Lambda$ for edges
incident to $V$ is contained in a vertex group of $\Lambda''$
which is contained in $\Lambda(V)$.  A generator of $\Lambda(V)$ which is
not contained in any edge group incident to $V$ is an element of
only that vertex group of $\Lambda$, but is also an element of
some vertex group of $\Lambda''$ and that vertex group of
$\Lambda''$ can only be contained in $\Lambda(V)$.  Thus we will
get that each generator of $\Lambda(V)$ is in a vertex group of
$\Lambda''$ which is contained in $\Lambda(V)$.

Consider then an edge group $C$ of an edge incident to $V$ in $\Lambda$. 
Delete the edges $E_1, \ldots E_n$ of the underlying tree of $\Lambda$ 
that are incident to $V$ with the edge group $C$, 
leaving a connected component $T_0$ containing $V$, and connected components
$T_1,\ldots,T_n$ with $T_i$ containing the vertex $V_i$ of $E_i$ opposite $V$ for each $i$. 
Then $W=A*_CB$ where $A$ is the group generated by the vertex groups of
the tree $T_0$ and $B$ is the group generated by the vertex groups of 
the forest $T_1,\ldots,T_n$. 
Neither $A$ nor $B$ equals $C$ as $\Lambda$ is reduced.  Each vertex group of
$\Lambda''$ is contained in a unique vertex group of $\Lambda$ and
so is contained in either $A$ or $B$ but not in both, since the
intersection of $A$ and $B$ is the edge group $C$. 
There is at least one vertex group of $\Lambda''$ in each of $A$ and $B$.
Hence there are adjacent vertices of $\Lambda''$ having vertex
groups one in $A$ and one in $B$, whose intersection is the edge group 
of $\Lambda''$ for the edge between these vertices.  
But the intersection of these vertex groups of
$\Lambda''$ is also contained in $C$.  Since the edge groups of
$\Lambda$ and $\Lambda''$ are visual subgroups conjugate to
$\langle J\rangle$, we have that $C$ is the edge group of $\Lambda''$ for the
edge between these vertices of $\Lambda''$. Hence $C$ is contained
in a vertex group of $\Lambda''$ contained in $A$.  
If this vertex group of $\Lambda''$ is contained in a vertex group in $A$ 
other than $V$, then $C$ would be contained in the edge groups in a
geodesic path between between $V$ and this other vertex in $T_0$.  
But the edge groups for edges incident to $V$ in $T_0$ are different
conjugates of $\langle J\rangle$ than $C$, since we deleted all edges
incident to $V$ having $C$ as edge group. 
As $C$ cannot be contained in a different conjugate of $\langle J\rangle$, 
instead the vertex group of $\Lambda''$ that is in $A$ and contains $C$ is actually
contained in the vertex group $\Lambda(V)$ of $V$ in $\Lambda$.
This completes the analysis of the claim and so completes the
proof of the lemma. 
\end{proof}

We now finish the proof of Theorem \ref{Lambdadecomptheorem}. 
By Lemma \ref{LambdaLambdapplemma}, the vertex groups of $\Lambda''$
are in fact equal to the vertex groups of $\Lambda$.
Finally we compare $\Lambda$ and $\Lambda''$ with $\Lambda'$. 
Each vertex group $G$ of $\Lambda''$ is a subgroup 
of a conjugate of a vertex group $G'$ of
$\Lambda'$ which is in turn a subgroup of a conjugate of a vertex
group $H$ of $\Lambda$.  But $G$ is a vertex group of $\Lambda$
and cannot be contained in a conjugate of another vertex group of
$\Lambda$ (since again $\Lambda$ is reduced).  Hence $G=H$,
$G^*\preceq {G'}^*\preceq G^*$ so $G^*={G'}^*$, and each vertex
group of $\Lambda$ is conjugate to a vertex group of $\Lambda'$.
On the other hand, if $H'$ is a vertex group of $\Lambda'$ then
$H'$ is a subgroup of a conjugate of a vertex group $H$ of $\Lambda$.  
But $H$ is also a vertex group of $\Lambda''$, is
contained in a conjugate of a vertex group $G'$ of $\Lambda'$, so
$G'=H'$ is conjugate to $H$. Hence the vertex groups of $\Lambda'$
correspond to conjugate vertex groups of $\Lambda$, as required for (3). 
This completes the proof of Theorem 8.1. 
\end{proof}

\begin{lemma} 
Suppose $\Lambda$ is a reduced visual graph of groups
decomposition for a Coxeter system $(W,S)$ of finite rank such
that the edge groups of $\Lambda$ are conjugates.  Then $\Lambda$
can be twisted resulting in a new generating set $S_*$ for $W$ and
a reduced visual graph of groups decomposition $\Psi$ such that
the edge groups of $\Psi$ are all equal.
\end{lemma}

\begin{proof}
If not, take an example of a $\Lambda$ for $(W,S)$ and an edge group $E_1$, 
having a minimum number of edges labelled by groups different from $E_1$, 
which cannot be twisted to a $\Psi$ with equal edge groups. 
Then some vertex group $V$ contains incident edge groups $E_1$ and $E_2$ with $E_1\neq E_2$. 
Let $W_1$ be the group generated by the vertex groups of $\Lambda$ that are joined to $V$ 
by a geodesic path in the underlying tree of $\Lambda$ terminating in an edge 
incident to $V$ labelled by $E_1$ 
and let $W_2$ be the group generated by the rest of the vertex groups of $\Lambda$. 
Then we have a free product decomposition $W = W_1*_{E_1}W_2$ 
with $E_2\subset V \subset W_2$. 
Let $d \in W$ be of minimal length, with respect to $S$, such that $E_1 =dE_2d^{-1}$. 
By considering the normal form for $d$ with respect to the 
amalgamated product $W_1*_{E_1}W_2$, we deduce that $d\in W_2$. 

Suppose $E_1=\langle S_0\rangle$,
$E_2=\langle \bar S_0\rangle$, $W_1=\langle S_1\rangle$, and
$W_2=\langle S_2\rangle$ for $S_0$, $\bar S_0$, $S_1$, $S_2\subset S$. 
By Lemma 4.3, we have $d\bar S_0d^{-1}=S_0$. 
Then twist the visual decomposition
$\langle S_1\rangle*_{\langle S_0\rangle}\langle S_2\rangle$ 
by conjugating the generators $S_2$ by $d$,
giving $S_*=S_1\cup dS_2d^{-1}$, and corresponding $\Lambda_*$
where we conjugate each vertex and edge group of $\Lambda$ with generators in $S_2$. 
An edge labelled $E_1$ cannot be in the twisted part of $\Lambda$
since the generators $S_0$ would have to be contained in each
vertex and edge group in a geodesic path between such an edge and an edge
with label $E_1$ incident at $V$.  Hence all the edges having
label $E_1$ originally still have label $E_1$ in $\Lambda_*$.  
The edge incident to $V$ labelled $E_2$ in the original $\Lambda$ is conjugated by $d$
to $dE_2d^{-1}=E_1$ and so we have at least one more edge labelled by
$E_1$, and hence at least one fewer edge labelled by a group
different from $E_1$.  Thus $\Lambda_*$ contradicts the minimality
of number of edges labelled by groups different from a particular
edge group assumed for $\Lambda$.  Instead, twisting to reduce the
number of edges labelled by a group different from a chosen edge
group must eventually transform a given $\Lambda$ to a graph of
groups $\Psi$ having all the same edge groups.
\end{proof}

\begin{lemma}  
Suppose $\Psi$ is a reduced graph of groups decomposition for a
Coxeter system $(W,S)$ of finite rank such that all of the
edge groups of $\Psi$ are equal. Suppose $\Psi'_0$ is a similar
decomposition for $(W,S')$ such that each vertex group is
conjugate to a vertex group of $\Psi$ and the equal edge groups of
$\Psi_0'$ are conjugate to the edge groups of $\Psi$. Then by a
sequence of twists applied to $\Psi'_0$ there results a new set of
generators $S'_*$ and corresponding visual graph of groups $\Psi'$
such that the vertex groups of $\Psi'$ are equal to those of
$\Psi$ and the edge groups of $\Psi'$ are all equal and equal to
the edge groups of $\Psi$, and hence $\tilde\Psi=\tilde\Psi'$.
\end{lemma}

\begin{proof}
Let $\tilde T$ be the Bass-Serre tree for $\tilde\Psi$.  
Then each vertex group $V'$ of $\Psi'_0$ stabilizes a v-vertex of $\tilde T$,
but stabilizes at most one v-vertex since $V'$ cannot be a subgroup
of a conjugate of an edge group of $\Psi'_0$, and the same is true for any
$\Psi'_*$ resulting by twists conjugating vertex groups and
preserving the same edge groups from $\Psi'_0$. Let $T_0$ be the
spanning tree for the v-vertices of $\tilde T$ that are stabilized by a vertex
group of such a $\Psi'_*$ and take $\Psi'_*$ so that $T_0$ has a
minimal number of vertices.  The smallest $T_0$ can be is one
v-vertex for each vertex group of $\Psi'_*$ plus one e-vertex, 
corresponding to the common edge group of $\tilde\Psi$, 
connected to each of the v-vertices of $T_0$.  
In this case, conjugating $\Psi'_*$ carries $T_0$ to the standard transversal $T_1$ 
of $\tilde T$ and so takes $\Psi'_*$ to a $\Psi'$ having the same vertex 
and edge groups as $\Psi$.

Suppose instead that $T_0$ has more than one e-vertex.  
Suppose further that some v-vertex $wV$ of $T_0$, for $V$ a vertex group of $\Psi$,
stabilized by a vertex group $V'=wVw^{-1}$ of $\Psi'_*$ has more
than one edge of $T_0$ incident at that vertex.  Let $E$ be the
common edge group of $\Psi$ so there are e-vertices $uE$ and $vE$
adjacent to $wV$,  the edge group $E'$ of $\Psi'_*$ is
$uEu^{-1}=vEv^{-1}$, and $uv^{-1}\in V'$. Twist $\Psi'_*$ to
$\Psi'_{**}$ by conjugating each vertex group of $\Psi'_*$
stabilizing a v-vertex of $T_0$ on the $vE$ side of $wV$ by the
element $uv^{-1}$.  Then $uv^{-1}E'vu^{-1}=uEu^{-1}=E'$, and so edge
groups have not changed.  

If $V'_2$ is a vertex group of $\Psi'_\ast$ stabilizing a
v-vertex $w_2V_2$ on the $vE$ side of $wV$, then $uv^{-1}V'_2vu^{-1}$
stabilizes $uv^{-1}w_2V_2$. If $p$ is a geodesic path from $wV$ to
$w_2V_2$ in $T_0$, then translating $p$ by $uv^{-1}$ results in a
path from $uv^{-1}wV=wV$ to $uv^{-1}w_2V_2$.  Since the first edge
in $p$ is $vE$, the first edge in the translated path is $uE$.  
We conclude that the spanning tree for the v-vertices stabilized by
$\Psi'_{**}$ consists of the part of $T_0$ on the $vE$ side of
$wV$ translated by $uv^{-1}$ together with the rest of $T_0$.
Since the e-vertex $vE$ is carried to the e-vertex $uE$ in the new spanning tree, 
there are fewer vertices in the new spanning tree, contradicting the
minimality of $T_0$ for $\Psi'_*$. 

Finally suppose that all the v-vertices of $T_0$ that are stabilized by a vertex
group of $\Psi'_*$ are the leaves of $T_0$ (the end points of $T_0$), and that
$T_0$ has at least two e-vertices.  Then $T_0$ has 
a v-vertex that is not stabilized by a vertex group of $\Psi'_*$.  
Let $\tilde T'$ be the Bass-Serre tree of $\tilde\Psi'_*$. 
Let $\tilde T'_*$ be the result of replacing, equivariantly with respect to the 
action of $W$, each translate of the standard transversal $T_1'$ 
of $\tilde T'$ by a copy of $T_0$
so that $T_0$ is attached by identifying each vertex of $T_0$ stabilized by a
vertex group $V'$ of $\Psi'_*$ with the vertex $V'$ of $T_1'$ (this vertex  
remains labelled $V'$). 
In particular, the e-vertices of $\tilde T'$ (those 
that are labelled by cosets of the edge group of $\Psi'_*$) are replaced 
by copies of the level one core of $T_0$ (the tree $T_0$ minus its leaves 
and their adjoining edges).  
Then $W$ acts on the tree $\tilde T'_*$ translating the vertices labelled by cosets of vertex
groups of $\Psi'_*$ in the same way as in $\tilde\Psi'_*$.  

Define a map $\tau:\tilde T'_* \to \tilde T$ by mapping  
the vertex $V'$ of the attached $T_0$ to the vertex $wV$ of $T_0$ 
when a vertex group $V'$ of $\Psi'_*$ stabilizes the vertex 
$wV$ in $T_0$, and by mapping the translates of $T_0$ in $\tilde
T'_*$ isomorphically to corresponding translates of $T_0$ in $\tilde T$ so as to
make $\tau$ respect the action of $W$. 
Then $\tau$ is locally injective,  
since the cosets $uv'E'$ of the edge group $E'$ in a given v-vertex $uV'$ of $\tilde T'$, 
which is also a vertex of $\tilde T'_*$,  
correspond to the cosets $uwvE$ of the edge group $E$ in the v-vertex 
$\tau(uV') = uwV$ of $\tilde T$ under the correspondence $v'=wvw^{-1}$, 
since $E' = wEw^{-1}$ and $V' = wVw^{-1}$. 
Hence $\tau$ is injective, since $\tau$ is a map of trees.  
But $T_0$ has an interior v-vertex $U = tV$, corresponding to a
vertex group $V$ of $\Psi$, which is not stabilized by a vertex group of $\Psi'_*$.  
Now $U = \tau(U')$ where $U'$ is an interior v-vertex of 
the attached $T_0$ in $\tilde T'_*$. 
Let $w^{-1}V'$ be the v-vertex of $\tilde T'$ stabilized by $V$. 
Then $w^{-1}V'w = V$, and so $V'=wVw^{-1}$. 
Hence the vertex group $V'$ of $\Psi'_*$ stabilizes 
the vertex $wV$ in $T_0$, and so $\tau(tw^{-1}V')=tw^{-1}wV=tV=U$. 
As $tw^{-1}V' \neq U'$ in $\tilde T'_*$, 
we have a contradiction to $\tau$ being injective. 
\end{proof}

Applying these lemmas to the result of the last theorem we have
the following conclusion.

\begin{theorem} 
{\rm (The Decomposition Matching Theorem)} 
Suppose $(W,S)$ and $(W,S')$ are Coxeter systems for the same
finitely generated Coxeter group and $W$ has a nontrivial
splitting as $A*_CB$. Then there are sequences of twists applied
to $(W,S)$ and $(W,S')$ giving rise to Coxeter systems $(W,S_*)$
and $(W,S'_*)$, respectively, such that there exists a nontrivial
reduced visual graph of groups decomposition $\Psi$ of $(W,S_*)$
and a nontrivial reduced visual graph of groups decomposition
$\Psi'$ of $(W,S'_*)$ having the same graphs and the same vertex
and edge groups and all edge groups equal and a subgroup of a
conjugate of $C$.
\end{theorem}

\begin{proof}
Take $\Lambda$ and $\Lambda'$ from Theorem \ref{Lambdadecomptheorem}, 
twist $\Lambda$ to get a visual decomposition $\Psi$ of $(W,S_*)$ 
with one edge group and twist $\Lambda'$ to get a visual decomposition 
$\Psi'_0$ with one edge group.  
Then twist $\Psi'_0$ to a visual decomposition $\Psi'_1$ of $(W,S'_*)$ having
the same vertex and edge groups as $\Psi$, so $\tilde\Psi=\tilde\Psi'$.  
Now $\Psi$ and $\Psi'_1$ only differ by expansions and contractions
rearranging the edge group attachments to vertex groups.  
So there is a visual decomposition $\Psi'$ of $(W,S'_*)$
with the same graph of groups structure as $\Psi$.
\end{proof}

\section{The Simplex Matching Theorem}  

The next lemma is known to experts. 
For a proof see Paris \cite{Paris}. 

\begin{lemma}
Let $W$ be a finitely generated Coxeter group 
with two complete Coxeter systems $(W,S)$ and $(W,S')$. 
Let 
$$(W,S) = (W_0,S_0)\times(W_1,S_1)\times\cdots\times(W_k,S_k)$$
with $(W_0,S_0)$ finite and $(W_i,S_i)$ infinite and irreducible for each $i=1,\ldots,k$. 
Suppose 
$$(W,S') = (W_0',S_0')\times(W_1',S_1')\times\cdots\times(W_\ell',S_\ell')$$
with $(W_0',S_0')$ finite and $(W_j',S_j')$ infinite and irreducible for each $j=1,\ldots,\ell$. 
Then $W_0 = W_0'$.  Let $Z = Z(W_0)$. 
Then $k = \ell$ and after reindexing we have 
$ZW_i = ZW_i'$ for each $i = 1,\ldots,k$. 
\end{lemma}

\begin{lemma}
Let $W$ be a finitely generated Coxeter group with two Coxeter systems 
$(W,S)$ and $(W,S')$. 
Let $S_1\subset S$ and $S_1'\subset S'$ 
and suppose that $W_1 = \langle S_1\rangle =\langle S_1'\rangle$. 
If the basic subgroups of $(W,S)$ isomorphically match basic subgroups of $(W,S')$, 
then the basic subgroups of $(W_1,S_1)$ isomorphically match basic subgroups of $(W_1,S_1')$.
\end{lemma}

\begin{proof} 
On the contrary, suppose $(\langle B_1\rangle,B_1)$ is a basic subgroup of $(W_1,S_1)$ 
that matches with a nonisomorphic basic subgroup $(\langle B_1'\rangle,B_1')$ of $(W_1,S_1')$. 
Without loss of generality, we may assume that $|\langle B_1\rangle|>|\langle B_1'\rangle|$. 
Then either $B_1$ is of type ${\bf C}_{2q+1}$ and $B_1'$ 
is of type ${\bf B}_{2q+1}$ for some $q\geq 1$ or 
$B_1$ is of type ${\bf D}_2(4q+2)$ and 
$B_1'$ is of type ${\bf D}_2(2q+1)$ for some $q\geq 1$. 
Let $B$ be the base of $(W,S)$ containing $B_1$ 
and let $B'$ be the base of $(W,S')$ matching $B$. 
Then $B$ is not of type ${\bf A}_5$ and $B_1$ 
is not of type ${\bf C}_2$. 
By Theorem 6.10, there is a $B_1''\subset B'$ such that $(\langle B_1''\rangle,B_1'')$ 
is a finite irreducible subgroup of $(W,S')$ 
and $[\langle B_1\rangle,\langle B_1\rangle]$ is conjugate to 
$[\langle B_1''\rangle,\langle B_1''\rangle]$ in $W$; 
moreover $(\langle B_1\rangle,B_1)\cong (\langle B_1''\rangle,B_1'')$, 
since $(\langle B\rangle,B)\cong (\langle B'\rangle,B')$. 
Now $B_1'$ is conjugate to $B_1''$ in $W$ by Theorem 6.12. 
Therefore $|\langle B_1\rangle| = |\langle B_1''\rangle| = |\langle B_1'\rangle|$ 
which is a contradiction. 
\end{proof}

\begin{theorem}
{\rm (The Simplex Matching Theorem)}
Let $(W,S)$ and $(W,S')$ be finite Coxeter systems 
with isomorphic matching basic subgroups. 
Then $(W,S)$ and $(W,S')$ have the same number of visual 
subgroups of each complete system isomorphism type. 
In particular, $|S| = |S'|$. 
\end{theorem}

\begin{proof}
The proof is by induction on $|S|$. 
This is clear if $|S| = 1$, 
so assume $|S| > 1$ and the theorem is true for all Coxeter systems 
with fewer generators than $|S|$. 
Assume first that $(W,S)$ is complete. 
Then $(W,S')$ is complete by Prop. 5.10. 
Let 
$$(W,S) = (W_1,S_1)\times\cdots\times(W_n,S_n)$$
be the factorization of $(W,S)$ into irreducible factors and suppose 
$(W_i,S_i)$ is finite if and only if $i\leq k$. 
Let 
$$(W,S') = (W_1',S_1')\times\cdots\times(W_m',S_m')$$
be the factorization of $(W,S')$ into irreducible factors and suppose 
$(W_i,S_i')$ is finite if and only if $i\leq \ell$. 
By Lemma 9.1, we have
$$(W_1,S_1)\times\cdots\times(W_k,S_k) = (W_1',S_1')\times\cdots\times(W_\ell',S_\ell').$$
By the Matching Theorem for systems of a finite Coxeter group, 
we can reindex so that $W_i$ is noncyclic if and only if $i\leq p$ 
and $W_j'$ is noncyclic if and only if $j\leq p$ 
and $[W_i,W_i] = [W_i',W_i']$ for each $i\leq p$. 
By hypothesis, $(W_i,S_i) \cong (W_i',S_i')$ for each $i\leq p$. 
As the remaining finite factors have order 2, we have $k = \ell$ 
and $(W_i,S_i) \cong (W_i,S_i')$ for $p < i\leq k$. 

By quotienting out the finite normal subgroup $(W_1,S_1)\times\cdots\times(W_k,S_k)$, 
we may assume that $W_i$ and $W_j'$ are infinite for each $i$ and $j$. 
By Lemma 9.1, we have that $m =n$ and after reindexing $W_i = W_i'$ for each $i$. 
Hence we may assume that $W$ is infinite and $(W,S)$ and $(W,S')$ are irreducible. 
By Theorem 2.4, we have that $(W,S) \cong (W,S')$. 
Thus in general $(W,S) \cong (W,S')$ when $(W,S)$ is complete. 

Now assume $(W,S)$ is incomplete. 
Then there are $a, b$ in $S$ such that $m(a,b) = \infty$. 
Hence 
$$W = \langle S-\{a\}\rangle \ast_{\langle S-\{a,b\}\rangle}\langle S-\{b\}\rangle$$
is a nontrivial visual amalgamated decomposition. 
By the Decomposition Matching Theorem, Theorem 8.6, 
there exist four nontrivial reduced visual graph of group decomposition of $W$, 
visual with respect to different sets of generators, a $\Lambda$ 
with respect to $S$, a $\Lambda'$ with respect to $S'$, 
a $\Psi$ with respect to another set of Coxeter generators $R$ of $W$, 
and a $\Psi'$ with respect to another set of Coxeter generators $R'$ of $W$ 
such that (1) the edge groups of $\Lambda$ and $\Lambda'$ are all conjugate 
and conjugate to a subgroup of $\langle S-\{a,b\}\rangle$; 
(2) there is a 1-1 correspondence between the vertices of $\Lambda$ and the vertices 
of $\Lambda'$ such that each vertex group of $\Lambda$ is conjugate to the 
corresponding vertex group of $\Lambda'$; 
(3) $\Psi$ is a twisted form of $\Lambda$ having all edge groups equal 
and conjugate to the edge groups of $\Lambda$, 
and having vertices in a 1-1 correspondence with those of $\Lambda$ 
such that each vertex group of $\Psi$ is conjugate to the corresponding vertex group 
of $\Lambda$, and $\Psi'$ is similarly a twisted form of $\Lambda'$; 
(4) $\Psi'$ is the same graph of groups as $\Psi$ and differs from $\Psi$ 
only in being a visual graph of groups decomposition of $W$ with respect 
to a different set of Coxeter generators. 

The Coxeter systems $(W,R)$ and $(W,S)$ are twist equivalent 
and so have the same number of visual subgroups of each complete 
system isomorphism type. 
Moreover $(W,R)$ and $(W,S)$ have isomorphic matching basic subgroups. 
Likewise the Coxeter systems $(W,R')$ and $(W,S')$ 
have the same number of visual subgroups of each complete 
system isomorphism type, and $(W,R')$ and $(W,S')$ 
have isomorphic matching basic subgroups. 

Let $\{(W_i,R_i)\}_{i=1}^k$ be the Coxeter systems of the vertex groups of $\Psi$, 
and let $(W_0,R_0)$ be the Coxeter system of the edge group of $\Psi$. 
Then $k\geq 2$, $R = \cup_{i=1}^kR_i$, and $\cap_{i=1}^k R_i= R_0$, 
and $R_i-R_0 \neq \emptyset$ for each $i > 0$, and 
$m(a,b) = \infty$ for each $a$ in $R_i-R_0$ and $b$ in $R_j-R_0$ with $i\neq j$. 
Let $\{(W_i',R_i')\}_{i=1}^k$ be the Coxeter systems of the vertex groups of $\Psi'$ 
indexed so that $W_i' = W_i$ for each $i$, 
and let $(W_0,R_0')$ be the Coxeter system of the edge group of $\Psi'$. 
Then $W_0'=W_0$, $R'= \cup_{i=1}^kR_i'$, and $\cap_{i=1}^k R_i'= R_0'$, 
and  $R_i'-R_0 \neq \emptyset$ for each $i > 0$, and 
$m(a',b') = \infty$ for each $a'$ in $R_i'-R_0'$ and $b'$ in $R_j'-R_0'$ with $i\neq j$. 
Moreover $(W_i,R_i)$ and $(W_i, R_i')$ have 
isomorphic matching basic subgroups for each $i$ by Lemma 9.2. 

Let ${\cal C}$ be a complete system isomorphism type and let ${\cal C}(S)$ 
be the number of visual subgroups of $(W,S)$ of isomorphism type ${\cal C}$. 
By the induction hypothesis, ${\cal C}(R_i) = {\cal C}(R_i')$ for each $i$. 
Observe that 
\begin{eqnarray*}
{\cal C}(S) & = & {\cal C}(R) \\
& = & \sum_{i=1}^k{\cal C}(R_i)-(k-1){\cal C}(R_0) \\
& = & \sum_{i=1}^k{\cal C}(R_i')-(k-1){\cal C}(R_0') \\
& = & {\cal C}(R') \ \ = \ \ {\cal C}(S'),
\end{eqnarray*}
which completes the induction.
\end{proof}

\section{The Maximum Rank of a Coxeter Group}

In this section we describe an algorithm for 
constructing a Coxeter system of maximum rank for 
a finitely generated Coxeter group.  
Let $(W,S)$ be a Coxeter system of finite rank. 
We say that $(W,S)$ can be {\it blown up along a base} $B$ if 
$(W,S)$ and $B$ satisfy the hypothesis of either Theorem 5.6 or 5.9. 
If $(W,S)$ can be blown up along a base $B$, 
then we can blow up $(W,S)$ to a Coxeter system $(W,S')$ as in 
the statement of Theorem 5.6 or 5.9 such that $|S'| = |S|+1$, 
the base $B$ matches a base $B'$ of $(W,S')$ 
with $|\langle B\rangle| >  |\langle B'\rangle|$,  
and each other base $C$ of $(W,S)$ is also a base of $(W,S')$. 
We say that $(W,S')$ is obtained by {\it blowing up} $(W,S)$ {\it along the base} $B$. 

By the process of blowing up along a base, 
we can effectively construct a sequence 
$S = S^{(0)}, S^{(1)},\ldots, S^{(\ell)}$ of Coxeter generators of $W$ 
such that $(W,S^{(i+1)})$ is obtained by blowing up $(W, S^{(i)})$ 
along a base for each $i = 0,\ldots, \ell-1$ and $(W,S^{(\ell)})$ cannot be blown up along a base. 
The sequence terminates since the sum of the orders 
of the basic subgroups decreases at each step of the sequence.  
By the next theorem, the system $(W, S^{(\ell)})$ has maximum rank 
over all Coxeter systems for $W$.

\begin{theorem}
{\rm (The Maximum Rank Theorem)}
Let $(W,S)$ be a Coxeter system of finite rank. 
Then the following are equivalent:
\begin{enumerate}
\item We have $|S| \geq |S'|$ for every set of Coxeter generators $S'$ of $W$. 
\item Each base $B$ of $(W,S)$ matches a base $B'$ of $(W,S')$ 
with $|\langle B\rangle| \leq |\langle B'\rangle|$ for every 
set of Coxeter generators $S'$ of $W$. 
\item The system $(W,S)$ cannot be blown up along a base. 
\end{enumerate}
\end{theorem}

\begin{proof}
Suppose that $|S| \geq |S'|$ for every set of Coxeter generators $S'$ of $W$ 
and on the contrary, a base $B$ of $(W,S)$ matches a base $B'$ of $(W,S')$ 
with $|\langle B\rangle| > |\langle B'\rangle|$. 
By the Basic Matching Theorem either 
$B$ is of type ${\bf C}_{2q+1}$ and $B'$ 
is of type ${\bf B}_{2q+1}$ for some $q\geq 1$ or
$B$ is of type ${\bf D}_2(4q+2)$ and 
$B'$ is of type ${\bf D}_2(2q+1)$ for some $q\geq 1$. 
By Theorems 5.7 and 5.8, we have that $(W,S)$ and $B$ 
satisfy the hypothesis of Theorem 5.6 or 5.9. 
Therefore $(W,S)$ can be blown up along $B$ to obtain 
a system $(W,S')$ with $|S'| = |S|+1$ contrary to the maximality of $|S|$. 
Therefore (1) implies (2). 

Suppose that each base $B$ of $(W,S)$ matches a base $B'$ of $(W,S')$ 
with $|\langle B\rangle| \leq |\langle B'\rangle|$ for every 
set of Coxeter generators $S'$ of $W$. 
If $|\langle B\rangle| = |\langle B'\rangle|$ for every base $B$ of $(W,S)$, 
then $(W,S)$ and $(W,S')$ have isomorphic matching basic subgroups by 
the Basic Matching Theorem, and so  
$|S| = |S'|$ by the Simplex Matching Theorem 9.3. 

Suppose a base $B$ of $(W,S)$ matches a base $B'$ of $(W,S')$ 
with $|\langle B\rangle| < |\langle B'\rangle|$. 
By the Basic Matching Theorem either 
$B'$ is of type ${\bf C}_{2q+1}$ and $B$ 
is of type ${\bf B}_{2q+1}$ for some $q\geq 1$ or
$B'$ is of type ${\bf D}_2(4q+2)$ and 
$B$ is of type ${\bf D}_2(2q+1)$ for some $q\geq 1$. 
By Theorems 5.7 and 5.8, we have that $(W,S')$ and $B'$ 
satisfy the hypothesis of Theorem 5.6 or 5.9, 
and so $(W,S')$ can be blown up along $B'$ 
to obtain a system $(W,S'')$ with $|S''| = |S'| +1$ 
such that $B'$ matches a base $B''$ of $(W,S'')$ 
with $|\langle B'\rangle| > |\langle B''\rangle| = |\langle B\rangle|$ 
and for each other base $C'$ of $(W,S')$ the base $C'$ is a base of $(W,S'')$. 

If $(W,S)$ and $(W,S'')$ do not have isomorphic matching basic subgroups, 
we can blow up $(W,S'')$ along a base. 
Continuing in this way, we obtain a sequence of Coxeter generators $S' = S^{(1)},\ldots, S^{(\ell)}$ 
of $W$ such that $(W,S^{(i+1)})$ is obtained by blowing up $(W, S^{(i)})$ 
along a base for each $i = 1,\ldots, \ell-1$ 
and $(W,S)$ and $(W,S^{(\ell)})$ have isomorphic matching basic subgroups. 
In particular, $|S^{(i+1)}| = |S^{(i)}|+1$ for each $i=1,\ldots,\ell-1$. 
By the Simplex Matching Theorem, $|S| = |S^{(\ell)}|$, and so $|S| > |S'|$. 
Thus (2) implies (1). 

Finally (2) and (3) are equivalent by the Basic Matching Theorem 
and Theorems 5.7 and 5.8. 
\end{proof}

We end our paper with the following theorem 
that says that any two Coxeter systems of maximum rank 
for a finitely generated Coxeter group have a lot in common. 

\begin{theorem}
{\rm (Simplex Matching Theorem for Maximum Rank Systems)}
Suppose $W$ is a finitely generated Coxeter group 
and $S$ and $S'$ are Coxeter generators of $W$ of maximum rank. 
Then $(W,S)$ and $(W,S')$ have the same number of visual 
subgroups of each complete system isomorphism type. 
\end{theorem}

\noindent\begin{proof} 
By the Maximum Rank Theorem, each base $B$ of $(W,S)$ 
matches a base $B'$ of $(W,S')$ with $|\langle B\rangle| =|\langle B'\rangle|$. 
Therefore $(W,S)$ and $(W,S')$ have isomorphic matching basic subgroups 
by the Basic Matching Theorem. 
Hence $(W,S)$ and $(W,S')$ have the same number of visual 
subgroups of each complete system isomorphism type 
by the Simplex Matching Theorem. 
\end{proof}

\end{document}